# Mathematical model for sustainable fisheries resource management accounting for size spectrum

Hidekazu Yoshioka[1, *], Yumi Yoshioka[2], Motoh Tsujimura[3], Ayumi Hashiguchi[4]

[1]Japan Advanced Institute of Science and Technology, 1-1 Asahidai, Nomi, Ishikawa 923-1292, Japan. ORCID: 0000-0002-5293-3246

[2] Gifu University, Yanagido 1-1, Gifu, Gifu 501-1193, Japan.

[3] Doshisha University, Karasuma-Higashi-iru, Imadegawa-dori, Kamigyo-ku, Kyoto, Kyoto 602-8580, Japan. ORCID: 0000-0001-6975-9304

[4] Okayama University, 3-1-1, Tsushima-Naka, Kita-ku, Okayama, Okayama 700-8530, Japan.

* Corresponding author: yoshih@jaist.ac.jp

***Abstract***

This paper proposes a novel modelling and control framework for growth models that incorporate a size spectrum in conjunction with numerical computation and extensive field surveys. In fisheries management, the size spectrum—characterized by individual differences in body weight and length—is a critical factor, as it influences the physiology and ecology of fish, as well as the preferences of anglers. However, a comprehensive theoretical framework for fisheries modelling and management that accounts for the size spectrum has yet to be established. We apply a growth model that considers the size spectrum to *Plecoglossus altivelis altivelis* (*Ayu*), an important inland fisheries resource in Japan. Additionally, we introduce a novel stochastic control theory for the resource management of *Ayu*, taking its size spectrum into account. The growth model is calibrated using data collected annually from a river system in Japan. Our control problem addresses the size spectrum of fishing benefits and terminal utility (nonlinear expectation) for sustainability, resulting in a nonstandard problem to which the dynamic programming principle does not apply. We address this difficulty using a time-inconsistent formalism, where solving the control problem is reduced to finding an appropriate solution to a system of nonlinear partial differential equations. We numerically compute the system using the finite difference method and explore the fisheries management of *Ayu* at the study site.

***Statements & Declarations***

**Acknowledgments:** We received assistance from the Hii River Fisheries Cooperative during field surveys. To conduct this research, we obtained special permission from the government of Shimane Prefecture (Tokusai No. 8 in 2025) to collect and capture the fish species *Plecoglossus altivelis altivelis* from the Hii River system in May and June 2025.

**Funding:** This study was supported by the Japan Science and Technology Agency (PRESTO No. JPMJPR24KE) and the Japan Society for the Promotion of Science (KAKENHI No. 25K00240).

**Conflict of Interests:** The authors declare no conflicts of interest.

**Data Availability:** The data will be made available upon reasonable request to the corresponding author.

**Declaration of Generative AI in Scientific Writing:** The authors did not use AI to generative the scientific results of this manuscript.

# 1. Introduction

## 1.1 Background

Biological resources, such as crops, fish, and livestock, are essential food sources for human life. Fish also play a crucial role in aquatic environments and ecosystems and have been frequently studied from an environmental perspective [1]. Enhancing our understanding of the status of fish as a biological resource is vital for implementing effective and sustainable management policies [2–4].

Fish living in natural environments have physiological spectra, such as individual differences in body weight and length [5,6], which are factors that shape various aspects of the food web, including mortality, fecundity, and feeding [7]. A small body size has been suggested to affect competition for space if the natural selection for a large body size is sufficiently strong [8]. Body size is a biological indicator of anthropogenic stressors [9]. The modality of body length distribution of pikeperch has been evaluated by considering environmental variables such as water temperature and precipitation [10]. Fish fecundity is considered an increasing function of body weight as implemented in operational models [11]. It has been suggested that catchability in fisheries may be affected by the body size [12]. The size spectrum also affects angler utility; it increases as the abundance or realization of larger fish increases [13,14].

Size spectrum of fish evolves over time because individual fish have different growth trajectories, which are considered to follow certain mathematical models, such as Von Bertalanffy [15,16], logistic or Ricker [17,18], hybrid models [19,20], or more sophisticated models considering the energy budget [21]. The size spectrum has been incorporated into a growth curve by randomizing its parameter values [22] or by adding a noise factor [23], which can be formulated through solutions to ordinary differential equations (ODEs) with random coefficients [24–26], and more recently, using an efficient formula to track the time-dependent size spectrum [27]. These methodologies have been employed in diverse research areas, such as population ecology [28], chemics [29], and hydrodynamics [30]; however, their applications in fisheries science are still limited to a few species, such as Atlantic herring [31] and *Plecoglossus altivelis altivelis* (commonly known as *Ayu* in Japan) [32,33], although both are important regional fishery resources, which is the motivation for this study.

Another motivation for this study is that adaptive (i.e., feedback or closed-loop) control theory for sustainable fisheries management has not been fully explored. Conventional optimal control theories are either deterministic or stochastic and are based on the total biomass or population. Solving a problem often reduces to finding a proper solution to an optimality equation [34–38]. These issues do not account for the size spectrum of fish or angler preferences. Moreover, formulating a control problem considering these elements results in nonstandard (known as "time-inconsistent" in control theory) ones as suggested in the literature on mathematical economics [39], followed by analytically tractable cases [40–42]. This is because of an increase in the problem dimension due to the existence of size spectrum that does not appear in classical control problems. From an application perspective, time inconsistency implies that dynamic programming (or the maximum principle) [43] serving as a standard tool in optimal control theory is no longer valid unless a special mathematical structure is assumed [44]. Therefore, a different methodology is required.

A possible way to resolve the aforementioned difficulty regarding the failure of dynamic programming is to use the notion of equilibrium control (Chapters 15 and 16 in Björk et al. [45]) to find a control to maximize the target objective function, which reduces the control problem to a system of nonlinear partial differential equations (PDEs), known as an extended Hamilton–Jacobi–Bellman (HJB) system. Technically, finding an equilibrium control requires studying the variational condition to be satisfied along with an extended HJB system. Because of the versatility of equilibrium control, this approach has a wider applicability range than the conventional control theory and has been applied to linear-quadratic games [46], portfolio management with utility-dependent discounts [47], mean–variance self-protecting insurance [41], insurance under heterogeneous time preferences [48], and environmental management with present bias [49]. Alderborn and Marín-Solano [50] formulated a long-term management problem for fisheries resource using multiple agents with nonexponential discounts. At present, we have not found an application study on equilibrium control in fisheries science that accounts for size spectrum of individual fish and angler preferences.

### 1.2 Aim and contribution

The aims of this study are as follows:

(i) To apply a growth curve model with a size spectrum to *Ayu* as a major fishery resource in Japan.

(ii) To formulate a time-inconsistent stochastic control theory for sustainable fisheries management using a size spectrum.

The contributions to achieving aims (i) and (ii) are explained below. See also **Table 1** summarizing the proposed model.

**(i) Growth curve model**

Our growth curve model is based on Yoshioka et al. [33] who proposed a tractable model for the body weight of *Ayu* in a river system based on a Von Bertalanffy model, where the asymptotic body weight (or equivalently, the maximum body weight) is assumed to be random to represent the size spectrum among individuals. We generalize this model so that other growth curve models, such as logistic ones, can be considered. Our generalization does not lose the analytical tractability of the original model that the probability density function (PDF) and its moments are determined explicitly. This property is useful for parameter estimations using real data.

Another innovation in this study from the previous study [33] is that the latest data collected in 2025 by the authors covers a longer period of body weight data of *Ayu*. This fish species has a one-year life history (**Section 4**), and previous studies covered their body weight data only from summer to autumn, which corresponds to the middle to late growth stages. By contrast, the data in this study cover the spring, corresponding to earlier growth stages, which potentially allows us to identify the parameter values of the growth curve model with higher accuracy. We also successfully collected a multiyear dataset of intensive surveys each summer, with which we could estimate the PDF of the growth curve model. Moreover, we obtained the length–weight relationship for inland fish [51–53] of the target fish at our study site such that the size spectra of body weight and length were used interchangeably.

Other growth curve models may also be used, such as diffusion process models [54–57] in which individuals grow with noisy trajectories. These models would be candidates for describing size spectrum if growth trajectories are observable with a sufficiently fine temporal resolution, which would be daily or finer; however, implementing such an observation scheme in a natural environment is technically difficult, and moreover, our target fish species *Ayu* often dies once it is captured. Our data are thus not chronological and that different data points represent different individuals [58]. We believe that temporally fine observation schemes be implementable in tailored artificial environments such as laboratories.

**(ii) Time-inconsistent control theory**

Based on the time-inconsistent control theory [45], we consider the control problem of a representative angler as a decision maker (this is not a social planner). We assume that he/she harvests *Ayu* during the harvesting period to maximize the total harvested biomass (cumulative benefit), but they may account for the sustainability of the fish such that a part of the fish remains at the end of the harvesting period (terminal utility). Both the cumulative benefit and terminal utility depend on the growth curve model and, hence, on the size spectrum.

Related studies that account for the size spectrum of fish are those on aquaculture systems, where the size spectrum is used to analyze cost-effective harvesting policies [59–61]. The main difference between the present and their formulations is that the former explores a closed-form harvesting policy with environmental uncertainties, and the timing of harvesting is adaptive, whereas the latter determines the harvesting amount at fixed time instances. Another difference between them is that the present study focuses on fishery resources in natural environments, whereas those in the literature focus on artificial environments.

Our control formalism, particularly that of terminal utility, is based on Desmettre and Steffensen [39] where preference heterogeneity among agents is provided by a nonlinear expectation that cannot be handled by dynamic programming, which is the reason we need the time-inconsistent control theory. Their methodology has been elaborated in insurance and finance [62,63], but not in other areas including resource management problems. We demonstrate that the equilibrium control for harvesting *Ayu* is determined by solving an extended HJB system, which is addressed computationally in this study using a finite difference method that preserves certain uniform bounds about its solution by selecting a sufficiently small time increment. With this numerical method, we finally apply the control theory along with the identified growth curve model for fisheries management of *Ayu*. Notably, for a limited problem class, such as the mean–variance terminal condition case, the extended HJB system can be solved by considering a continuum of time-consistent (and hence standard) HJB equations [64,65]; however, this methodology does not seem to apply to our problem.

### 1.3 Structure of this paper

**Section 2** presents and analyzes the growth curve model. **Section 3** presents and analyzes the control theory and derives the equilibrium control to maximize utility in fisheries management. **Section 4** presents the application of the growth curve model and control theory. **Section 5** summarizes the results obtained in this study and presents the perspectives of the study. **Appendix** presents the auxiliary results (**Appendix A**), numerical discretization for the extended HJB system (**Appendix B**), and technical results (**Appendix C**).

**Table 1.** Overview of the proposed model.

| Input | Growth model, objective |
|---|---|
| Objective | Profit by harvesting, terminal utility |
| Solution method | Extended HJB system |
| Advantage | Size spectrum can be considered |
| Disadvantage | More complex than the classical (time-consistent) control theory |

## 2. Fish growth model

### 2.1 Model description

Time is denoted as $t \geq 0$, where the initial time is 0. The body weight of an individual fish at time $t$ is $W_t$. Our growth curve model is given by

$$W_t = Kf(t), \quad t \geq 0, \tag{1}$$

where $K > 0$ is the asymptotic body weight that is assumed to be random and $f$ is a deterministic continuous function of time that is bounded between 0 and 1. Yoshioka et al. [33] used the Von Bertalanffy model (we scaled the parameters from the original model such that they are consistent with other models):

$$f(t) = \left(1 - \left(1 - (f_0)^{\frac{1}{3}}\right) e^{-\frac{r}{3}t}\right)^3, \quad t \geq 0 \tag{2}$$

with the growth rate $r > 0$ and initial condition $f_0 \in (0,1)$. This $f$ solves the ODE

$$\frac{\mathrm{d}f(t)}{\mathrm{d}t} = r\left(f(t)\right)^{\frac{2}{3}}\left(1 - \left(f(t)\right)^{\frac{1}{3}}\right), \quad t > 0 \tag{3}$$

with $f(0) = f_0$. The initial body weight is $W_0 = Kf_0$. The following logistic model is also used later:

$$f(t) = \frac{1}{\left((f_0)^{-1} - 1\right)e^{-rt} + 1}, \quad t \geq 0. \tag{4}$$

This $f$ solves the ODE

$$\frac{\mathrm{d}f(t)}{\mathrm{d}t} = rf(t)\left(1 - f(t)\right), \quad t > 0 \tag{5}$$

with $f(0) = f_0$. These models are common in (deterministic) growth curve descriptions because they can reproduce the sigmoid (i.e., S-shaped) growth trajectories of biological organisms. They share common parameters $f_0$ and $r$, and their total number of parameters is equal. We will also consider an extended model where $r$ is time-dependent in the ODE (5), that is, $r_t = r_0 + r_1 t$ with $r_0, r_1 \geq 0$.

### 2.2 Statistical properties

Yoshioka et al. [33] assumed that the asymptotic body weight $K$ is different among individuals such that the size spectrum is reproduced. Let the PDFs of $K$ and $W_t$ be $p_K$ and $p_{W_t}$, respectively. An elementary calculation yields

$$p_{W_t}(w) = \frac{1}{f(t)} p_K\left(\frac{w}{f(t)}\right), \quad w > 0. \tag{6}$$

Therefore, the PDF $p_{W_t}$ (and statistics) of body weight is available in closed form if $p_K$ is. We assume the following gamma distribution because it was successfully applied to the body weight data of *Ayu* until 2024 [33]:

$$p_K(k) = \frac{1}{\Gamma(\alpha)\beta^\alpha} k^{\alpha-1} e^{-\frac{k}{\beta}}, \quad k > 0, \tag{7}$$

where the shape parameter $\alpha > 0$ modulates the shape of the PDF, and the scale parameter $\beta > 0$ stretches it. The gamma distribution provides a unimodal PDF of the body weight used in our case study, and its statistics, such as moments, can be obtained explicitly (**Section 4**). We have

$$\int_0^{+\infty} w p_{W_t}(w) \mathrm{d}w = \alpha\beta f(t) \leq \alpha\beta f(T) = \int_0^{+\infty} w p_{W_T}(w) \mathrm{d}w, \quad 0 \leq t \leq T, \tag{8}$$

stating that the average body weight increases over time.

Our formulation suggests that the statistics of body length can also be obtained if we use the allometric relationship between body weight $w > 0$ and body length $l > 0$ [51]:

$$w = al^b \tag{9}$$

with parameters $a > 0$ and $b > 0$. The parameter $b$ is often close to 3 [50,52], reflecting the scaling relationship that body weight $w$ is approximately proportional to body volume $O(l^3)$. With this formula, we will check whether our data follow this scaling relationship.

## 3. Control theory

We present and analyze a resource management model.

### 3.1 System dynamics

We consider a simple control problem in which an angler, a union member of regional fisheries cooperatives, dynamically optimizes his arrival intensity at a target river. We consider a control problem in a harvesting period during which the fish do not have reproduction, and hence their population does not increase. We assume that their arrival at the river follows a counting process $P=(P_t)_{t\geq 0}$ and that their jump intensity, the arrival intensity $u=(u_t)_{t\geq 0}$, is adaptively controlled so that an objective function is maximized.

The state variable is the population $X=(X_t)_{t\geq 0}$ of fish in part of the river before reproduction. The population dynamics are assumed to be [66]

$$\underbrace{\mathrm{d}X_t}_{\text{Increment}} = \underbrace{-h(X_{t-})\mathrm{d}P_t}_{\text{Harvesting}} \underbrace{-\kappa X_{t-}\mathrm{d}D_t}_{\text{Catastrophe}},\ t>0, \tag{10}$$

starting from the initial condition, that is, $X_0>0$. In the first term on the right-hand side of (10) for harvesting events, $h$ is a nonnegative and Lipschitz continuous function with $0\leq h(x)\leq x$ for $x\geq 0$. We use $h(x)=\min\{\bar{h},x\}$ with a constant $\bar{h}>0$ with which the angler harvests fish for a constant amount $\bar{h}$ at each arrival except when the population is scarce ( $x<\bar{h}$ ), so that the population remains nonnegative after each jump of $P$. In the second term on the right-hand side of (10) for catastrophic events, $\kappa\in(0,1]$ is a positive constant and $D=(D_t)_{t\geq 0}$ is a counting process, where $P$ and $D$ are obtained via random time-changes (or thinning) of two independent standard Poisson random measures $N_P(\mathrm{d}t,\mathrm{d}z)$ and $N_D(\mathrm{d}t,\mathrm{d}z)$, respectively. The jump intensity $v=(v_t)_{t\geq 0}$ of $D$ is given by

$$v_t = d + ku_t^{\gamma},\ t>0 \tag{11}$$

with $d,k>0$ and $\gamma>1$.

At each jump of $D$, the population decreases as $X_t=(1-\kappa)X_{t-}$ and becomes extinct if $\kappa=1$ [66] and $\kappa\in(0,1)$ corresponds to a partial extinction. We assume that the rate of (partial) extinction is triggered by excessive predation pressures from anglers, which was modeled using a mean-field [66]. We do not consider the mean-field case because it is extremely complex and a time-inconsistent control case without mean-field effects is already more difficult than conventional stochastic control problems. The term $d$ in (11) represents catastrophe by natural origins, such as predation from waterfowls, and the other $ku_t^{\gamma}$ by harvesting, where the assumption $\gamma>1$ implies that the angler is averse to extinction such that selecting an extremely large $u_t$ is effectively avoided.

### 3.2 Control problem

Let $\mathbb{U}$ be the set of admissible controls $u=(u_t)_{t\geq 0}$:

$$\mathbb{U}=\left\{u=(u_t)_{t\geq 0}: 0\leq u\leq \bar{U},\ u \text{ is adapted to } \mathbb{F}, \text{ the SDE (10) admits a unique pathwise solution}\right\}, \tag{12}$$

where $\mathbb{F}$ is a natural filtration generated by $P$ and $D$, and $\bar{U}>0$ is a constant. From a practical perspective, $\bar{U}$ would be at most 1 (1/day) by assuming that an angler does not arrive at the river more than once per day. To simplify the notations, we set

$$\bar{W}_t \equiv \int_0^{+\infty} w p_{W_t}(w)\mathrm{d}w. \tag{13}$$

In the rest of **Section 3**, we assume that the harvesting fish is opened within the fixed time interval $(0,T)$ with $T>0$. We set the objective function $J_T$ to be maximized by a representative angler (the expectation $\mathbb{E}$ acts for $X$):

$$J_T = \underbrace{\mathbb{E}\left[\sum_{\substack{k=1,2,3,\ldots \\ 0<\tau_k<T}} \left(\int_0^{+\infty} w p_{W_{\tau_k}}(w)\,\mathrm{d}w\right) h\left(X_{\tau_k -}\right)\right]}_{\text{Profit by harvesting}} + \eta \underbrace{\int_0^{+\infty} p_{W_T}(w)\,\rho^{(-1)}\left(\mathbb{E}\left[\rho\left(wX_T\right)\right]\right)\mathrm{d}w}_{\text{Terminal utility}}. \tag{14}$$

Here, $\eta>0$ is a multiplier for balancing the two terms and $\rho:[0,+\infty)\to\mathbb{R}$ with $\rho(0)<\infty$ is an invertible function that is concave or convex.

The first term in the right-hand side of (14) is the profit by harvesting that is assumed to be proportional to the total harvested fish biomass, where the strictly increasing sequence $\{\tau_k\}_{k=1,2,3,\ldots}$ represents jump times of $P$. The second term in the right-hand side of (14) is the terminal utility as a certainty equivalent expectation (i.e., nonlinear expectation) version of $\int_0^{+\infty} p_{W_T}(w)\mathbb{E}[wX_T]\mathrm{d}w$, which evaluates the biomass of fish remaining at $T$. The number of eggs reproduced by some fish species, such as *Ayu*, at a single reproduction event is approximately proportional to its body weight [67], but with significant uncertainties (Fig.2(A) in Barneche et al. [68]), thereby motivating us to use a nonlinear expectation. A function $\rho$ used in this study is the power-type

$$\rho(y)=\frac{1}{\psi+1}y^{\psi+1} \text{ with } \rho^{(-1)}(y)=(\psi+1)^{\frac{1}{\psi+1}}y^{\frac{1}{\psi+1}},\ y\geq 0 \tag{15}$$

with shape parameter $\psi>-1$. Other convex $\rho$ are also possible, such as the exponential one [69,70]:

$$\rho(y)=\exp(\psi y) \text{ with } \rho^{(-1)}(y)=\frac{1}{\psi}\ln y,\ y>0 \tag{16}$$

with $\psi\neq 0$, but we do not focus on this case because (15) is more convenient for applications. With the exponential one (16), we must handle a huge $\psi y$, which becomes computationally intractable. For the power case, two body weight distributions having the same mean but different variances, skewness or tails, apparently would generate the same equilibrium control. However, this is not the case in general when $\rho$ is not of a power type. In our application study, we assume a power type for convenience, but theoretically, a wider range of $\rho$ can be considered. Note that even under the power case the corresponding extended HJB system is not analytically solvable to the best of the authors' knowledge. Moreover, solutions to the extended HJB system are highly nontrivial in this case.

With (15) or (16), we have

$$\rho^{(-1)}\left(\mathbb{E}\left[\rho\left(wX_T\right)\right]\right)\geq\mathbb{E}\left[wX_T\right]\ (\text{resp., } \leq) \tag{17}$$

if $\psi>0$ with convex $\rho$ (resp., $\psi<0$ with concave $\rho$), representing an overestimation (resp., underestimation) of reproductive success. The equality in (17) holds if $\psi=0$. Therefore, the value of $\psi$ determines how the angler evaluates terminal utility; selecting a smaller $\psi$ implies that they are more pessimistic about the fish biomass remaining at the terminal time. The terminal time $T$ corresponds to the end of the fishing season (pre-spawning period), and from a resource economics standpoint, the terminal term represents the non-use or conservation value of the remaining spawning stock biomass. The concave specification reflects risk aversion or diminishing marginal utility regarding the stock size left for future generations.

The coexistence of a linear expectation (first term) and a nonlinear expectation (second term) in (14) hinders us from applying the dynamic programming principle to find $u \in \mathbb{U}$ that maximizes $J_T$. The approach based on time-inconsistent control resolves this issue and gives us a way to find an equilibrium $u$.

***Remark 1.*** We have $J_T < +\infty$ for any $u \in \mathbb{U}$ because

$$\begin{aligned} J_T &\le \bar{W}_T h(X_0) \mathbb{E}\left[ \sum_{\substack{k=1,2,3,\ldots \\ 0<\tau_k<T}} 1 \right] + \eta \int_0^{+\infty} p_{W_T}(w) \rho^{(-1)}\left( \mathbb{E}\left[ \rho(wX_0) \right] \right) \mathrm{d}w \\ &\le \bar{W}_T \left( h(X_0) \bar{U} T + \eta X_0 \right) \\ &< +\infty \end{aligned} . \tag{18}$$

### 3.3 Equilibrium control

The formulation in this subsection is based on Chapters 15 and 16 of Björk et al. [45] with an adaptation to our case. We set a dynamic version of $J_T$:

$$J(t,x;u) = \mathbb{E}^{t,x}\left[ \sum_{\substack{k=1,2,3,\ldots \\ t<\tau_k<T}} \bar{W}_{\tau_k} h\left( X_{\tau_k -} \right) \right] + \eta \int_0^{+\infty} p_{W_T}(w) \rho^{(-1)}\left( \mathbb{E}^{t,x}\left[ \rho(wX_T) \right] \right) \mathrm{d}w \tag{19}$$

for $0 \le t \le T$ and $x \ge 0$, where $\mathbb{E}^{t,x}$ is the expectation conditioned on $X_t = x$. We have $J_T = J(0, X_0; u)$.

We consider the following control problem that maximizes $J$:

$$\Phi(t,x) = \text{"}\max_{u \in \mathbb{U}}\text{"} J(t,x;u) . \tag{20}$$

Here, "max" in (20) is understood within controls such that the following **Definition 1** is satisfied (equilibrium control) (Definition 15.3 in Björk et al. [45]). In the sequel, a Markovian control $u \in \mathbb{U}$ is a control that depends only on the current state, that is, formally $u_t = u(t, X_t)$.

***Definition 1***

*Consider a Markovian control $\hat{u} \in \mathbb{U}$. For any arbitrary Markovian control $u \in \mathbb{U}$ and a real number $\delta \in (0, T-t)$, set $u_\delta \in \mathbb{U}$ as*

$$u_\delta(s,y)=\begin{cases}u(s,y) & (t\le s<t+\delta,\ y\ge 0)\\ \hat{u}(s,y) & (t+\delta\le s\le T,\ y\ge 0)\end{cases},\quad 0\le t<T. \tag{21}$$

*If*

$$\liminf_{\delta\to+0}\frac{J(t,x;\hat{u})-J(t,x;u_\delta)}{\delta}\ge 0 \tag{22}$$

*for all* $u\in\mathbb{U}$, $0\le t<T$, *and* $x\ge 0$, *then we call* $\hat{u}$ *an equilibrium control and* $J(\cdot,\cdot;\hat{u})$ *the value function.*

The condition (22) states that an equilibrium control outperforms any other small modification $u_\delta$ under the limit $\delta\to+0$, and hence has a variational flavor.

We conclude this section with our extended HJB system, from which an equilibrium control is formally derived. First, we have the following boundary condition at $x=0$:

$$\Phi(t,0)=0,\quad 0\le t<T \tag{23}$$

and the terminal condition at $t=T$:

$$\Phi(T,x)=\eta\bar{W}_T x,\quad x\ge 0. \tag{24}$$

Second, we set the partial differential operator $\mathbb{L}_\theta$ parameterized by $\theta\in[0,\bar{U}]$:

$$\mathbb{L}_\theta\phi(t,x)=\theta\left(\phi(t,x-h(x))-\phi(t,x)\right)+\left(d+k\theta^\gamma\right)\left(\phi(t,x-\kappa x)-\phi(t,x)\right) \tag{25}$$

for $0\le t\le T$ and $x\ge 0$, and generic sufficiently smooth functions $\phi$. We write $\partial_t\phi=\frac{\partial\phi}{\partial t}$.

Let $\lambda$ be the derivative of $\rho^{(-1)}$. Based on Section 3 of Desmettre and Steffensen [39], our extended HJB system is given as follows: for any $0\le t<T$ and $x>0$,

$$\partial_t\Phi(t,x)+\max_{0\le\theta\le\bar{U}}\left\{\begin{array}{l}\underbrace{\mathbb{L}_\theta\Phi(t,x)+\theta h(x)\bar{W}_t}_{\text{Coming from dynamics and harvesting}}\\ \underbrace{-\eta\left(\partial_t G(t,x)+\mathbb{L}_\theta G(t,x)\right)+\eta\int_0^{+\infty}p_{W_T}(w)\lambda\left(g(t,x,w)\right)\left(\partial_t g(t,x,w)+\mathbb{L}_\theta g(t,x,w)\right)\mathrm{d}w}_{\text{Coming from time-inconsistency of the terminal condition}}\end{array}\right\}=0, \tag{26}$$

where

$$G(t,x)=\int_0^{+\infty}p_{W_T}(w)\rho^{(-1)}\left(g(t,x,w)\right)\mathrm{d}w \tag{27}$$

and $g(\cdot,\cdot,w)$ for each $w>0$ solves

$$\partial_t g(t,x,w)+\mathbb{L}_{\hat{\theta}(t,x)}g(t,x,w)=0,\quad 0\le t<T \text{ and } x>0, \tag{28}$$

the boundary condition

$$g(t,0,w)=\rho(0),\quad 0\le t<T, \tag{29}$$

and the terminal condition

$$g(T,x,w)=\rho(wx),\ \ x\ge 0\ . \tag{30}$$

Here, for each $(t,x)$, we set

$$\hat{\theta}(t,x)=\underset{0\le\theta\le\bar{U}}{\arg\max}\left\{\begin{aligned}&\mathbb{L}_\theta\Phi(t,x)+\theta h(x)\bar{W}_t-\eta\left(\partial_t G(t,x)+\mathbb{L}_\theta G(t,x)\right)\\&+\eta\int_0^{+\infty}p_{W_T}(w)\lambda\left(g(t,x,w)\right)\left(\partial_t g(t,x,w)+\mathbb{L}_\theta g(t,x,w)\right)\mathrm{d}w\end{aligned}\right\}. \tag{31}$$

Our extended HJB system contains an HJB-like equation (26) and a continuum of backward Kolmogorov-like equations (28), both being interconnected through $\hat{\theta}$ in (31). Therefore, this system is more complicated than conventional HJB equations, which essentially have an equation of the form (26) with terms containing $g$ being dropped. This occurs when no terminal utility exists ($\eta=0$). The HJB-like equation (26) has terms coming from the nature of dynamic programming for harvesting benefit and those from the time inconsistency of the terminal utility that is based on the nonlinear expectation. These two elements appear separately inside "max" of (26) but are accounted for simultaneously in solving this equation.

According to **LemmaC1 in Appendix**, the added terms (second line) in the HJB-like equation (26) are nonnegative (resp., nonpositive) for convex (resp., concave) $\rho$, respectively if $g$ is non-negative. This implies how the nonlinear expectation enters the HJB-like equation. For convex $\rho$ with which the nonlinear expectation overestimates the corresponding (linear) expectation, the nonnegativity of the added terms possibly increases the quantity inside "max," and then formally yields a larger value function. A symmetric argument applies to concave $\rho$. The added terms with the auxiliary variables $g$ therefore correspond to the nonlinearity of the terminal condition in the extended HJB system.

***Remark 2.*** We formally have

$$\partial_t G(t,x)=\int_0^{+\infty}p_{W_T}(w)\lambda\left(g(t,x,w)\right)\partial_t g(t,x,w)\mathrm{d}w \tag{32}$$

by (27). Therefore, the two temporal derivative terms cancel in (26) and (31).

### 3.4 Verification result

We study the extended HJB system presented in the previous subsection, because it is the core of our control theory. We intend to show that equilibrium control in the sense of **Definition 1** is given by

$$\hat{u}_t=\hat{\theta}(t,X_t),\ \ 0<t<T\ . \tag{33}$$

The collection of all the functions that are continuous in $[0,T]\times[0,+\infty)$ and continuously differentiable for the first argument in $(0,T)$ is denoted by $C^{1,0}$. The verification results presented below demonstrate that control (33) based on the extended HJB system is an equilibrium control under certain conditions.

***Proposition 1***

*Suppose that the extended HJB system (26)-(31) admits a solution* $(\Phi, g)$ *such that* $\Phi \in C^{1,0}$ *and* $g(\cdot,\cdot,w) \in C^{1,0}$ *for all* $w>0$, $G \in C^{1,0}$, *and both* $\Phi$ *and* $g(\cdot,\cdot,w)$ *(*$w>0$*) grow at most polynomially for large* $x>0$. *Moreover, assume that*

$$\mathbb{E}^{t,x}\left[\left\{\Phi(t,X_t)\right\}^2\right],\ \mathbb{E}^{t,x}\left[\left\{g(t,X_t,w)\right\}^2\right] < +\infty,\ \textit{for any}\ t\in[0,T],\ x>0,\ w>0, \quad (34)$$

*and* $\hat{\theta} \in \mathbb{U}$. *Then,* $\hat{u}$ *in (33) is an equilibrium control in the sense of* ***Definition 1****.*

## 4. Application

### 4.1 Target fish

The target fish in this study are the diadromous fish *Ayu* that migrate between the Hii River system (a class-A river system in the eastern part of Shimane Prefecture, San-in region, Japan) and the Sea of Japan every year (**Figure 1**). The common life history of *Ayu* is completed in one year, which starts from hatching from eggs in a river during autumn, followed by wintering of larvae in shallow seas, juvenile upstream migration from the sea to a river, and growth in the midstream river environment in the river [71]. We focus on the life stage of *Ayu* in the river system (mid-April to the end of October) during which the fish are harvested in the Hii River system (July 1 to the end of October). Both natural and released *Ayu* exist in the Hii River system as in many other rivers in Japan, but their ratio is unknown.

The authors have collaborated with the Hii River Fisheries Cooperative (HRFC) to conduct field surveys on *Ayu*. From 2017 to 2024, a member of the HRFC collected daily fish catch and average body weight data from July 1 to late October to early November each year, which covered the growth data in summer and autumn that correspond to the growing to maturity periods, but not earlier data (May and June) that correspond to the juvenile period. In 2025, we started sampling individuals of *Ayu* from May with four members of the HRFC, such that earlier growth data could be collected being different from the previous years. We also conducted one-day intensive surveys (**Table 2**) around August each year, during which approximately 20 members of HRFC harvested *Ayu* and we measured individual body weights to study the size spectrum of the fish, which is $p_{W_t}$ in our context. We also obtained body length data for 2025, from which we estimated the allometric scaling relationship (9). In this section, the beginning of May 1 is set to the time 0 each year, around which *Ayu* migrate toward upstream along the river system.

We also conducted an almost weekly environmental DNA sampling analysis at fixed points near Kisuki and Shinmitoya in **Figure 1** from March to November 2025 (quantitative Polymerase Chain Reaction method [72], laboratorial analysis was due to Research Center of Environment & Pollution Co., Ltd.). Positive DNA concentrations were detected from April 17 (around the upstream migration) to November 7 (around the downstream migration) in 2025; the latter corresponded to approximately 190 days from May 1 in the aforementioned setting. This survey suggests that the growth dynamics of *Ayu* in the Hii River system should be tracked until the beginning of November, as we did in 2025.

We examine the Von Bertalanffy and logistic models. We identify the parameter values in the model for each year, which are $f_0 \in (0,1)$, $\alpha > 0$, $\beta > 0$, and $r > 0$, using the following strategy. Our strategy for parameter identification is based on that of Yoshioka et al. [33] with an improved treatment of sampled average weight data. We first identify the "shape" of the PDF $p_{W_t}$, which is essentially $p_K$ owing to (6), based on the body weight data at intensive surveys. Let $t_I$ be the time of intensive survey in some year, then we use the following moment matching method for average $\mathbb{E}\left[W_{t_I}\right]_{\text{Data}}$ and variance $\mathbb{V}\left[W_{t_I}\right]_{\text{Data}}$ of body weight data (in the sequel, subscripts "Data" and "Model" indicate quantities estimated using empirical data and the growth curve model, respectively):

$$\mathbb{E}\left[W_{t_I}\right]_{\text{Data}} = \mathbb{E}\left[W_{t_I}\right]_{\text{Model}} = \alpha\beta f(t_I) \text{ and } \mathbb{V}\left[W_{t_I}\right]_{\text{Data}} = \mathbb{V}\left[W_{t_I}\right]_{\text{Model}} = \alpha\beta^2 \left\{ f(t_I) \right\}^2 . \quad (35)$$

We can estimate $\alpha$ and $\beta$ from (35) as follows:

$$\alpha = \frac{\left(\mathbb{E}\left[W_{t_I}\right]_{\text{Data}}\right)^2}{\mathbb{V}\left[W_{t_I}\right]_{\text{Data}}} \text{ and } \beta = \frac{1}{f(t_I)} \frac{\mathbb{V}\left[W_{t_I}\right]_{\text{Data}}}{\mathbb{E}\left[W_{t_I}\right]_{\text{Data}}} . \quad (36)$$

Based on (36), we use the weighted least-squares method to fit $f_0$ and $r$:

$$\text{Minimize Err} = \frac{1}{\sum_{m=1}^{M} N_m} \sum_{m=1}^{M} N_m \left( \mathbb{E}\left[W_{t_m}\right]_{\text{Data}} - \mathbb{E}\left[W_{t_m}\right]_{\text{Model}} \right)^2 \text{ for } f_0 \in (0,1) \text{ and } r > 0 , \quad (37)$$

where the sequence $\{t_m\}_{m=1,2,3,\ldots,M}$ with $M \in \mathbb{N}$ represents the dates on which *Ayu* was caught in the Hii River system in the selected year. The weighting factor $N_m \in \mathbb{N}$ is the total number of samples caught on $t_m$. This weighting factor is coming from the estimate that the sum of $N_m$ independent and identically distributed random variables having average and variance approximates the true average with the error of $O\left(N_m^{-0.5}\right)$ (e.g., Example 5.21 in Klenke [73]. Intuitively, a data point with a larger $N_m$ would be more reliable in fitting the average. Having found fitted values of $f_0$ and $r$, we can find $\alpha$ and $\beta$ from (36).

***Remark 3.*** Because angling gear and mesh sizes may introduce size-selective sampling bias, the caught sample distribution may slightly differ from the true population distribution. This kind of uncertainties may be addressed by using the nonlinear expectation.

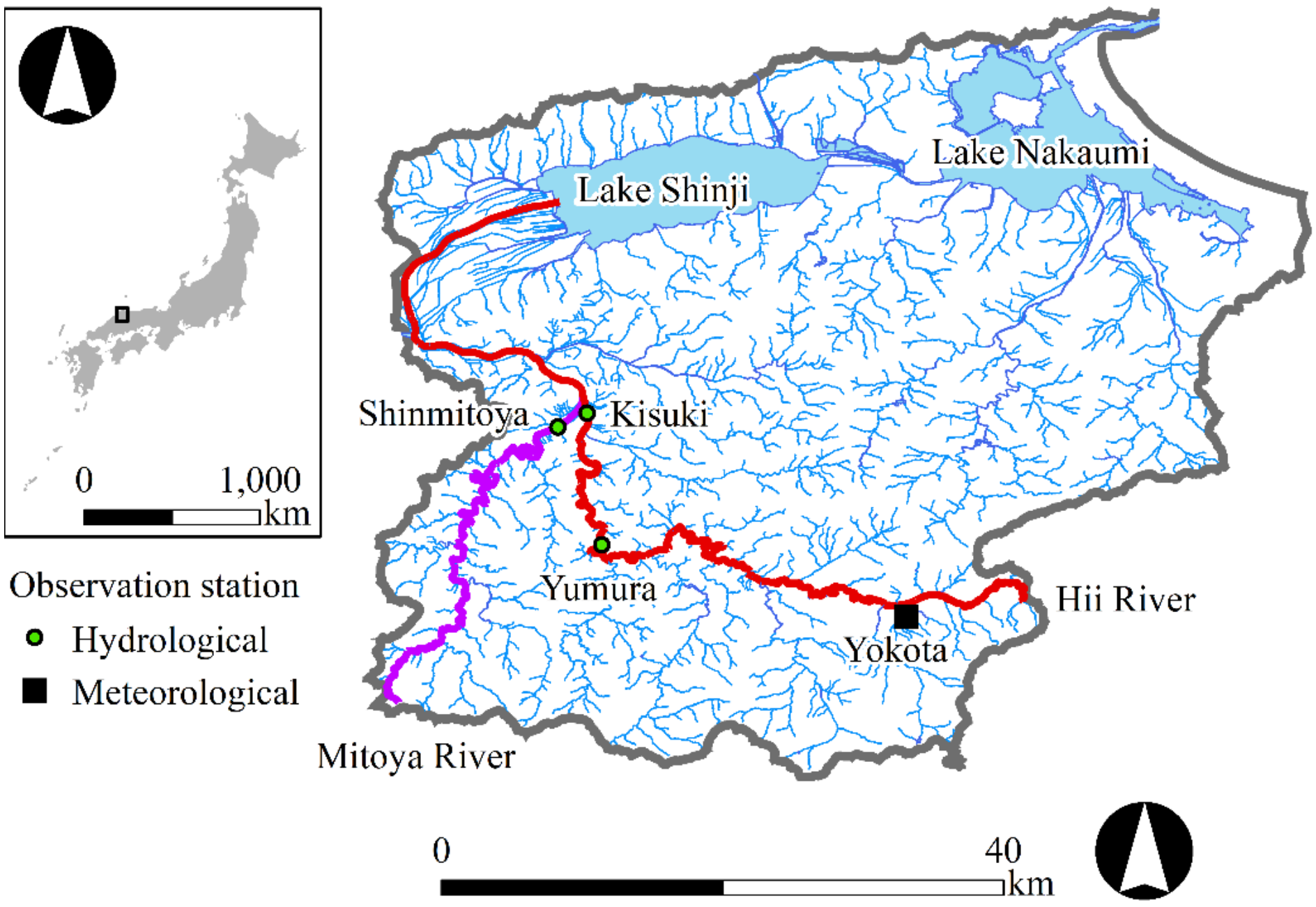

**Figure 1.** Map of the study site.

**Table 2.** Statistics and data on intensive surveys. The day 0 is May 1 each year; "na" indicates no data, "Std" denotes standard deviation, and "CV" denotes coefficient of variation. Intensive surveys were not conducted in 2020 through 2022.

| Year | 2016 | 2017 | 2018 | 2019 | 2023 | 2024 | 2025 |
|---|---|---|---|---|---|---|---|
| Date | Aug. 7 | Aug. 6 | Aug. 5 | Aug. 4 | Jul. 30 | Aug. 24 | Aug. 23 |
| Time (day) | 98 | 97 | 96 | 95 | 90 | 115 | 114 |
| No. of Samples | 207 | 234 | 189 | 227 | 297 | 459 | 446 |
| Average (g) | 55.2 | 55.6 | 57.3 | 56.4 | 52.2 | 49.9 | 48.2 |
| Std (g) | na | 19.1 | 18.5 | 18.2 | 21.0 | 15.7 | 16.7 |
| CV (-) | na | 0.34 | 0.32 | 0.32 | 0.40 | 0.31 | 0.35 |
| Skewness (-) | na | 0.77 | 1.15 | 0.95 | 1.41 | 1.52 | 1.44 |
| Median (g) | na | 52.8 | 54.5 | 54.0 | 46.5 | 47.0 | 44.2 |
| Maximum (g) | 120.5 | 132.0 | 152.0 | 119.5 | 163.0 | 127.4 | 139.0 |
| Minimum (g) | 38.0 | 20.5 | 16.0 | 20.0 | 11.0 | 13.3 | 14.7 |

## 4.2 Results and Discussion: parameter identification

### 4.2.1 Fitting results

**Figure 2** shows the empirical averages and theoretical averages of the fitted models (Von Bertalanffy and logistic) for each year. **Figure 3** shows the empirical and theoretical PDFs of the body weight during the intensive surveys. Both the Von Bertalanffy and logistic models provide an identical PDF in intensive surveys because of the parameter identification strategy explained in the previous subsection. **Tables 3 and 4** list fitted parameter values and the minimized least-squares error "MinErr" in (37) each year for Von Bertalanffy and logistic models, respectively.

### 4.2.2 Models before 2024

First, we discuss models before 2024. **Tables 3 and 4** suggest that the logistic model outperforms the Von Bertalanffy model for the years 2017–2019 and 2023–2024 considering the minimized least-squares error. Moreover, the theoretical value of the average body weight at the beginning of May 1, which is $\mathbb{E}[W_0]$ in these tables, is unrealistically small for the Von Bertalanffy model in 2018 and 2023, which is 0 (g) in 2023, and should not be possible. The logistic model for 2023 encounters the same issue but provides a more realistic result than that of the Von Bertalanffy model. The peculiar result about the theoretical initial body weight can be explained by form (2) such that taking $f_0 = 0$ does not cause any mathematical issues, whereas the logistic one (4) does. This difference is considered the reason for the more realistic initial body weight in the fitted logistic model. See also **Figure A2**.

For the growth rate $r$, its fitted values are in the order of $O(10^{-3})$ (1/day) to $O(10^{-2})$ or $O(10^{-1})$ (1/day) except for the years 2017 and 2024 in the Von Bertalanffy model, where $r$ is approximately $O(10^{-6})$ (1/day) to $O(10^{-5})$ (1/day), which is extremely small compared with those of other cases. According to **Figure A1**, the empirical average values are reasonably within or around the theoretical average ± standard deviation each year with slight overshoot of empirical values in 2023 and 2024. In summary, the comparison between the two models suggests that the logistic model is more accurate than the Von Bertalanffy model in our application study.

Finally, note that we proposed a parameter identification method based on the tractability of the model that does not require complex numerical computations. While the effectiveness of this method was suggested, it was also confirmed that unrealistic parameter values were estimated in some cases. Regarding this issue, it is desirable to conduct further examinations in the future, such as through application to other fish species.

### 4.2.3 Models in 2025

Using data from 2025, we examine the influence of the data because we could collect body weight data before opening the fishing of *Ayu* this year. **Figure 4(a)** compares the fitted model for three cases: with all the data discussed above, with data after the opening of fishing (after July 1), and with data before

September 1, which covers the intensive survey date. This figure suggests that a model identification procedure without using early-stage growth data (i.e., data from May 1 to June 31) significantly underestimates the growth rate, with a critical deviation between the theory and data at the early growth stage, whereas the other two cases provide almost the same average curves that are difficult to distinguish visually. Moreover, it also suggests that data after September are not considered useful for evaluating growth if those during the early growth state are available, as in 2025. This implies that data both before and after the inflection point of $f$ should be collected, if possible, for a better estimation of fish growth.

**Figure 2(f)** shows that the Von Bertalanffy and logistic models yield qualitatively different results, where the former predicts convex growth while the latter predicts sigmoidal growth. The empirical data imply that the logistic model is more reasonable by capturing the early growth stage, whereas the Von Bertalanffy model better minimizes the error, as summarized in **Tables 3 and 4**. We infer that the convex growth curves predicted by both models for 2024 are possibly owing to the lack of early growth data; nevertheless, verifying this hypothesis is impossible because such data do not exist.

We also examine another logistic model in which the growth rate $r$ is time-dependent in the ODE (5); i.e., $r_t = r_0 + r_1 t$ with $r_0, r_1 \geq 0$, which improves the accuracy of fitting the average body weight data. The fitting procedure is the same as that used in the previous case. The identified parameter values are summarized in **Table A1** for 2017–2019 and 2023–2025. Including the time dependence of the growth rate, both quantitatively and qualitatively, improves the accuracy of the average and standard deviation curves with better coverage of the empirical data (**Figures A1(f)** and **A2(f)**). **Figure 4(b)** shows a comparison of the fitted model under the data availability conditions, as explained above. **Figure 4(a)** suggests that ignoring the data after September 1 does not critically affect the identification results, whereas the lack of data during the early growth stage is crucial.

Finally, we identify the allometric relationship (9) using body weight (g) and length data (cm) collected in 2025 (500 samples), as shown in **Figure 5**, with the parameter values of $a$ and $b$ estimated using the common least-squares method. We have $a = 0.0054$ (g/cm$^b$) and $b = 3.19$ (-), suggesting that the identified $b$ is close to 3, as observed in a previous study of *Ayu* in a different river in Japan [52]. The collected data, identified models, and allometric relationships between body weight and length can be used as benchmarks for studying the yearly transitions of fish biology in the Hii River system. A multiyear dataset of body weight data of *Ayu*, similar to that collected in our study, is rare to the best of the authors' knowledge. Moreover, the ecology and biology (body weight and migration timing) of fish have been suggested to differ among different years in Japan [74], and the proposed model needs to be identified for each river (if data are available) at this stage. Establishing a normalization method to compare the growth of *Ayu* in different rivers is an important topic that should be addressed in the future.

#### 4.2.4 Relation to water temperature

**Table 2** suggests important findings regarding the body weight of *Ayu* in the Hii River system at intensive surveys. First, the empirical coefficient of variation is in the range of 0.3 to 0.4 and the empirical skewness at approximately 1 for all years. Moreover, **Figure 3** suggests that the PDF is unimodal for all years. These findings imply that the probabilistic nature of the body weight of *Ayu* in the Hii River system in summer has not critically changed in recent ten years. From the logistic model with the time-dependent growth rate, the average body weight at the beginning of October 31 is estimated as 98.3 (g), 90.3 (g), 101.7 (g), 77.8 (g), 66.9 (g), and 57.1 (g), from years 2017–2019 and 2023–2025, respectively.

For the relationship between water temperature and fish growth, Nakagawa et al. [75] reported that in the Kamo River in the Kansai Region in Japan, the water temperature exceeded 25 °C, which has been considered as the upper limit of suitable water temperature for *Ayu*, in downstream and midstream reaches. Rising water temperature would not only affect growth dynamics, particularly when it becomes outside the habitat suitability of target fish [76], but also the fish spawning phenology and consequently fisheries constraints [77] and food web dynamics [78]. Rising trends in water temperature have also been observed in other areas where commercially important fish species are migrating, such as North America [79] and Asia [80]. Based on the multimodal ensemble analysis, water temperature in global rivers has been predicted to increase by 1–5 °C from 2020 to 2099 [81]. According to these predictions, the rising trend of water temperature in the Hii River system may continue; however, the ranges are uncertain. Therefore, the fitted growth curve models and their spectra, particularly reproduction, which depends on body weight during autumn, could be considered subject to model uncertainties.

We consider that the trend toward miniaturization of *Ayu* in the Hii River system is due to rising water temperature, which is considered to lead to an increase in metabolic rate. The monthly air temperature rise in the study area (Yokota in **Figure 1**) was reported in Yoshioka et al. [33] from 2017 to 2024, and we reanalyzed it by adding data for 2025 (**Table A2 through A4**)[1]. As summarized in these tables, the air temperature was the first or second hottest in 2024 and 2025 from July to October, which corresponds to the growing season of *Ayu* in the Hii River system. Therefore, the rising trend in air temperature is suggested to continue until 2025 in the study area. In the Hii River system, we measured the river water temperatures at Yumura, Kisuki, and Shinmitoya in **Figure 1** from mid-May 2025 to the end of November 2025 (**Figure A3**). The obtained data reveal that the river water temperature in the midstream Hii River system exceeds 25 °C during summer and even exceeds 30 °C both in the mainstream and branch. Similar observations have been reported for Yumura in 2024 [33]. These observations are consistent with the hypothesis that high water temperature contributed to the lower body weight of *Ayu* in 2024 and 2025.

---

[1] Available at Japan Meteorological Agency https://www.data.jma.go.jp/stats/etrn/ (last accessed on November 7, 2025)

**Table 3.** Fitted parameter values and the minimized Err (MinErr): Von Bertalanffy model.

| | 2017 | 2018 | 2019 | 2023 | 2024 | 2025 |
|---|---|---|---|---|---|---|
| $\alpha$ (-) | 8.47.E+00 | 9.59.E+00 | 9.60.E+00 | 5.22.E+01 | 1.01.E+01 | 8.36.E+00 |
| $\beta$ (g) | 4.14.E+08 | 1.21.E+01 | 2.21.E+01 | 2.10.E+01 | 2.62.E+11 | 1.55.E+01 |
| $f_0$ (g) | 5.42.E-09 | 2.69.E-02 | 4.22.E-02 | 0.00.E+00 | 1.14.E-11 | 2.71.E-01 |
| $\mathbb{E}[W_0]$ (g) | 1.90.E+01 | 3.12.E+00 | 8.95.E+00 | 0.00.E+00 | 3.02.E+01 | 3.50.E+01 |
| $r$ (1/day) | 2.34.E-05 | 3.78.E-02 | 1.90.E-02 | 6.31.E-02 | 1.07.E-06 | 6.08.E-03 |
| MinErr ($g^2$) | 6.61.E+01 | 6.02.E+01 | 4.97.E+01 | 2.06.E+02 | 1.25.E+02 | 1.39.E+02 |

**Table 4.** Fitted parameter values and the minimized Err (MinErr): logistic model ($\alpha$ is the same with the Von Bertalanffy model).

| | 2017 | 2018 | 2019 | 2023 | 2024 | 2025 |
|---|---|---|---|---|---|---|
| $\alpha$ (-) | 8.47.E+00 | 9.59.E+00 | 9.60.E+00 | 6.18.E+00 | 1.01.E+01 | 8.36.E+00 |
| $\beta$ (g) | 1.56.E+02 | 1.03.E+01 | 1.32.E+01 | 1.28.E+01 | 1.02.E+05 | 5.76.E+00 |
| $f_0$ (g) | 1.73.E-02 | 1.14.E-01 | 1.07.E-01 | 3.46.E-02 | 3.08.E-05 | 6.53.E-02 |
| $\mathbb{E}[W_0]$ (g) | 2.28.E+01 | 1.12.E+01 | 1.37.E+01 | 2.75.E+00 | 3.18.E+01 | 3.15.E+00 |
| $r$ (1/day) | 9.47.E-03 | 2.48.E-02 | 1.99.E-02 | 4.42.E-02 | 3.91.E-03 | 1.12.E-01 |
| MinErr ($g^2$) | 6.48.E+01 | 5.91.E+01 | 4.72.E+01 | 1.97.E+02 | 1.25.E+02 | 1.43.E+02 |

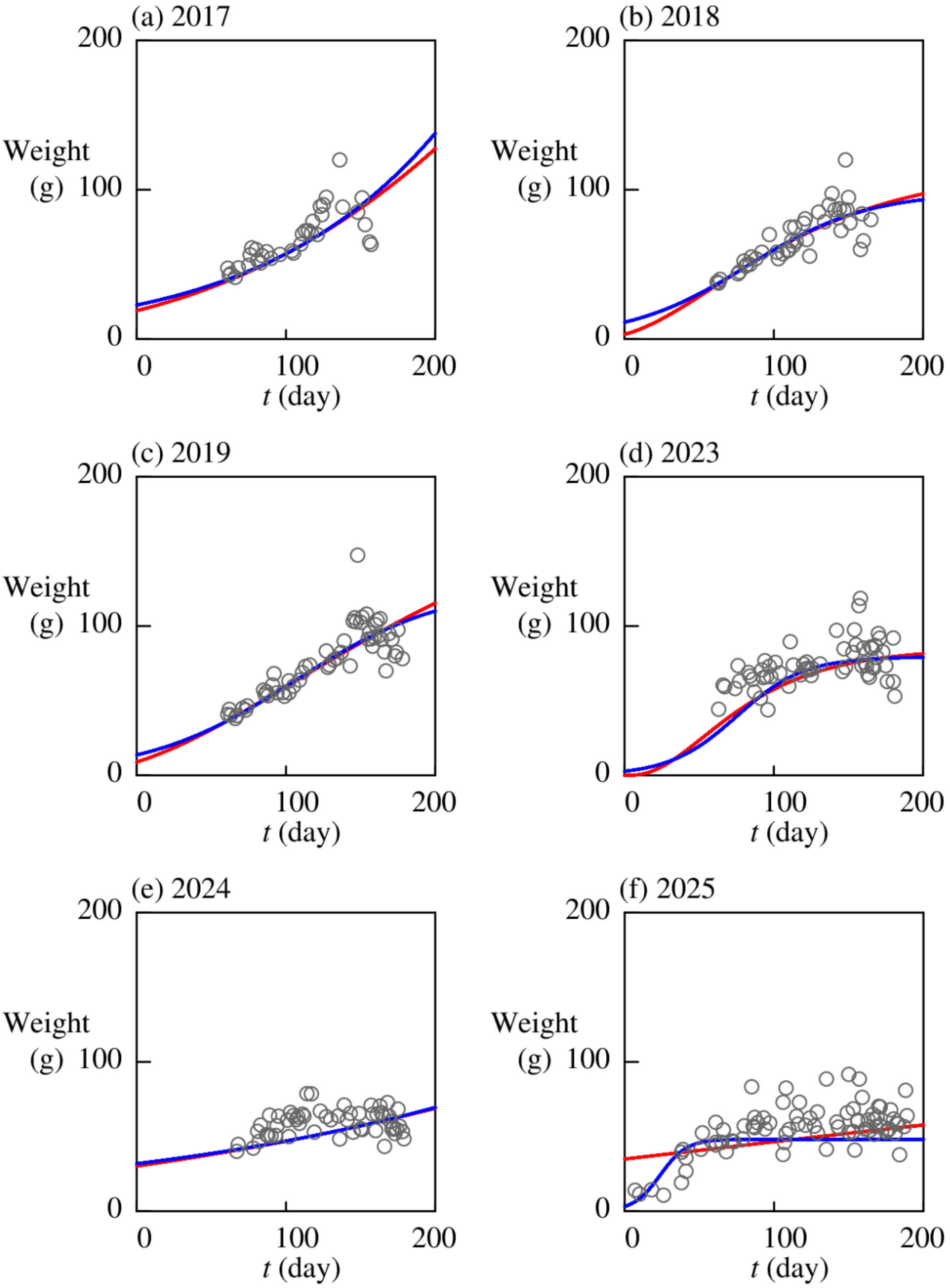

**Figure 2.** Comparison of average body weight each year: empirical (circles) and theoretical (curves) with Von Bertalanffy (red curves) and logistic models (blue curves).

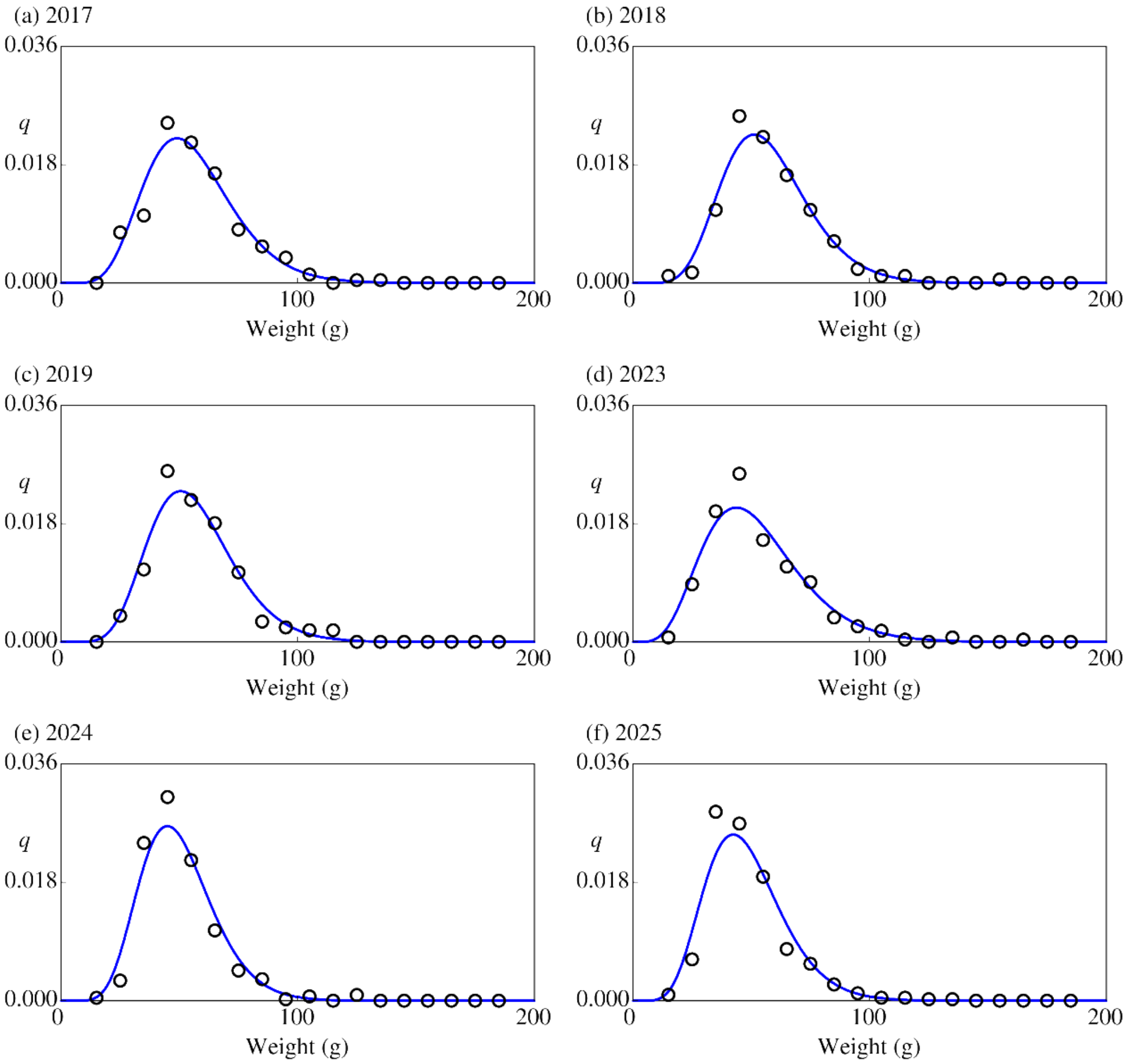

**Figure 3.** Comparison of the PDFs of body weight in intensive surveys each year: empirical (circles) and theoretical (curves).

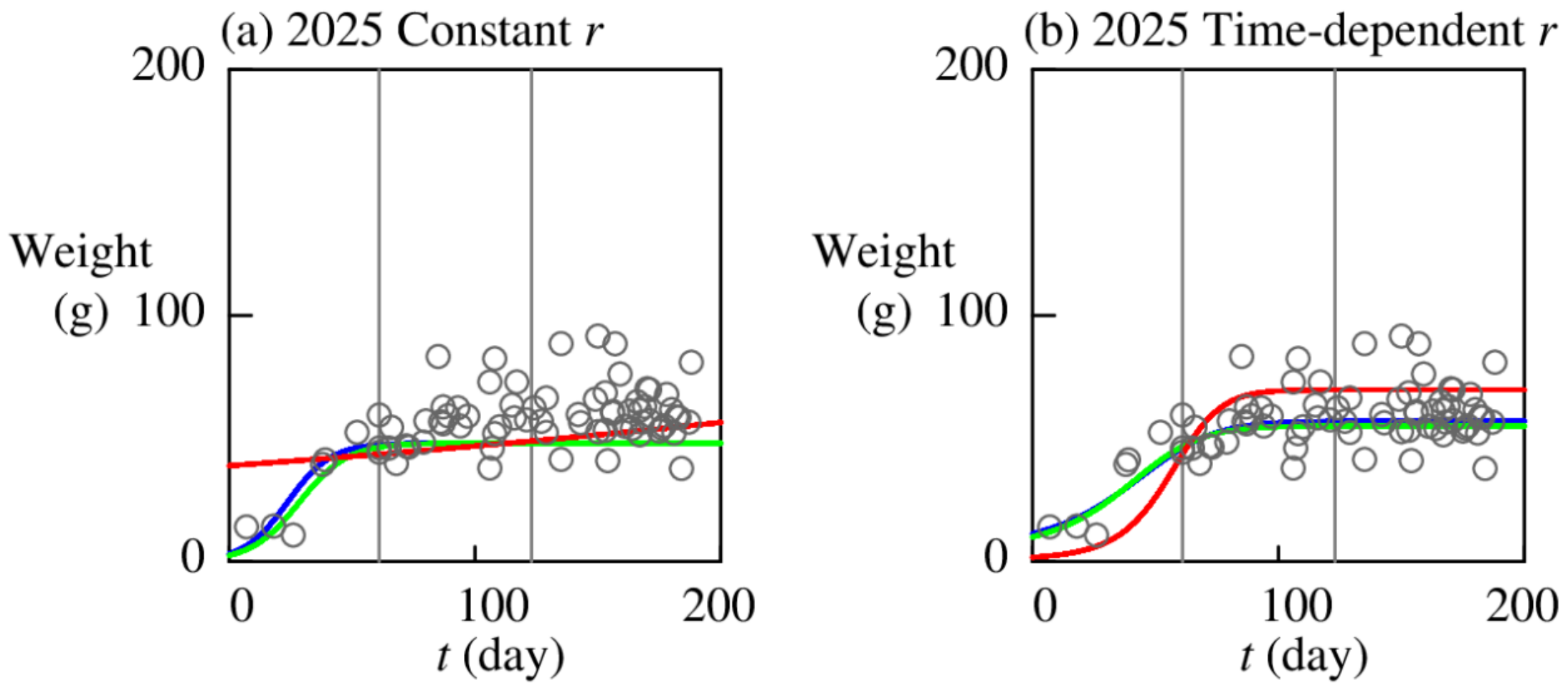

**Figure 4.** Comparison of average body weight each year with different fitting periods (a) logistic model with constant growth rate and (b) logistic model with variable growth rate: empirical (circles) and theoretical (curves) with fitted model using the entire data (blue curves), fitted model with data before September 1 (green curves), and fitted model with data after July 1 (red curves). Vertical gray lines represent the beginning of July 1 and September 1.

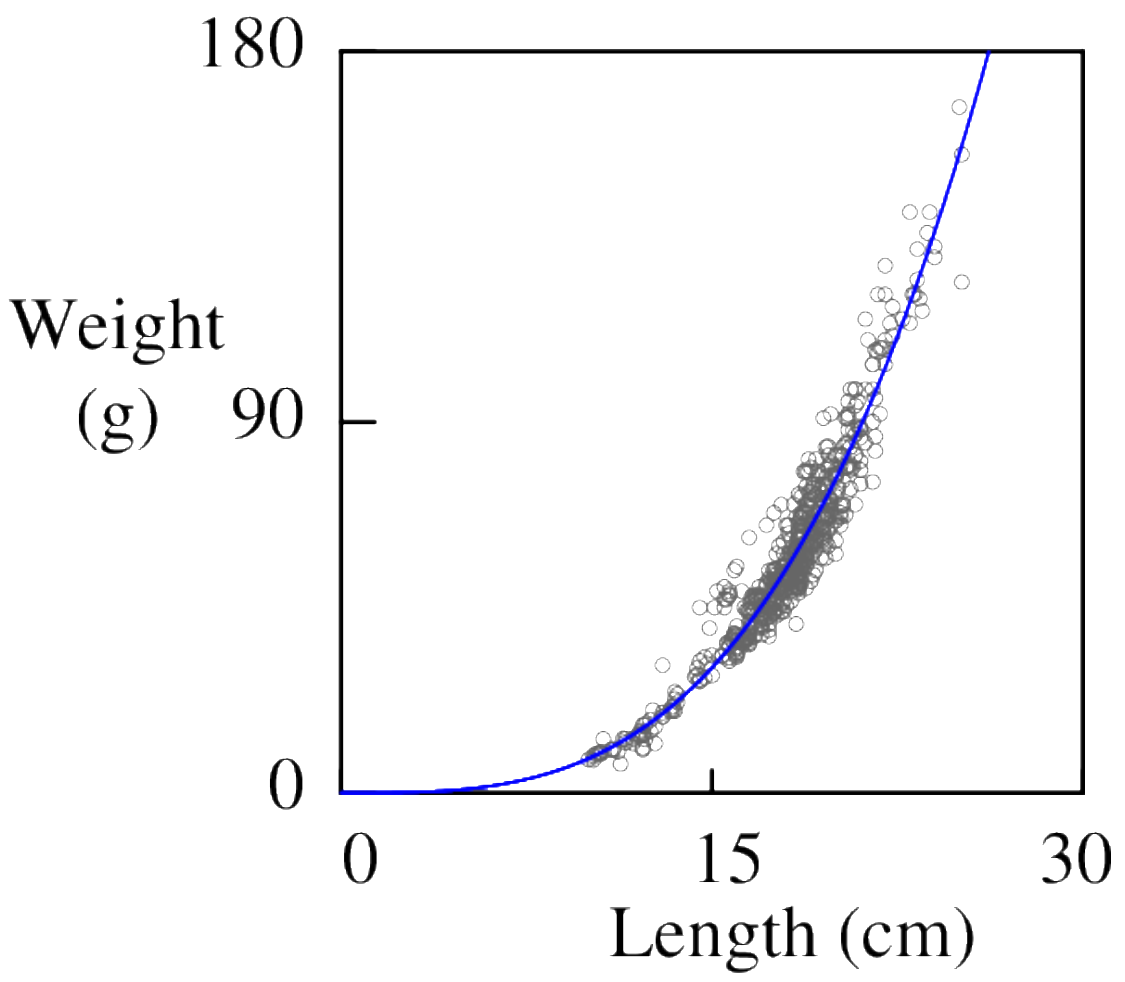

**Figure 5.** Allometric relationship between the body weight (g) and body length (cm) in 2025: empirical (circles) and theoretical (curve).

### 4.3 Application to fisheries management

We compute the equilibrium control based on the logistic model because of its superior accuracy compared with the Von Bertalanffy model, as reported in the previous subsection. The focus here is the dependence of the equilibrium control on the spectrum of body weight and anglers' preferences concerning the resource sustainability of *Ayu*.

We use the power-type $\rho$ in (15) and the following parameter values, considering the fishing history of a member in the HRFC who caught approximately 3,000 *Ayu* individuals each year: $\bar{U} = 1$ (1/day) (this constrains an average of one arrival per day and does not preclude multiple arrivals occurring on the same day), $\bar{X} = 4,000$, and $\bar{h} = 40$. We also set $d = 0.0001$ (1/day) assuming that the natural extinction probability is low, $k = 0.002$ (1/day) with $\gamma = 2$ assuming that the extinction rate due to harvesting dominates the natural extinction probability when the intensity of harvesting is high, and $\kappa = 1$ assuming that extinction is complete if it occurs. For the growth curve, we use a logistic model with a time dependent growth rate in 2025. The computational resolution is $\Delta t = 0.01$ (day) and $\Delta x = h$. The discretization in the $w$ direction is based on the quantile discretization suited to Gamma distributions (e.g., Section 2.3 in Yoshioka [82]). The terminal time is set as 181 (day) near the end of October and the 61 (day) corresponding to July 1. The specified terminal time is before the end of the downstream migration for spawning, according to the environmental DNA survey described in **Section 4.1**. We examine different values of the parameters $\eta$ and $\psi$ in the objective function. We checked sufficiency of the employed computational resolution. For example, doubling the current resolution in $w$ (199 grid points in this paper) to 399 and halving $\Delta t$ to 0.005 (day) does not critically change the value function; the relative change was 0.2% for the computational condition of **Figure 10(a)** and their difference was not visible.

**Figure 6** shows the computed value function $\Phi$ and the corresponding equilibrium control $\hat{u}$ when no terminal utility $\eta = 0$ is present (this case is independent of $\psi$). Similarly, **Figure 7** shows the computed $\Phi$ and $\hat{u}$ with $\eta = 0.6$ and $\psi = 0$. A comparison between **Figures 6 and 7** suggests that accounting for the terminal utility for the sustainability of the fisheries resource yields a less intensive arrival of the representative angler at the study site, which is visualized by the decrease in the red to blue areas in **Figure 7(b)** compared with those in **Figure 6(b)**. The numerical solutions are nonnegative, as theoretically suggested in **Proposition B1**, for both convex and concave $\rho$.

We investigate the influence of the parameters $\eta$ and $\psi$. **Figures 8 and 9** show the computed value functions $\Phi$ and $\hat{u}$, respectively, for different values of $\eta$. Similarly, **Figures 10 and 11** show the computed value functions $\Phi$ and control $\hat{u}$ for different values of $\psi$, which cover both the convex ($\psi > 0$) and concave ($\psi < 0$) $\rho$ cases. **Figure 8** suggests that the profile of the value function $\Phi$ becomes simpler as $\eta$ increases, that is, as the terminal utility for sustainability dominates the benefit of harvesting. **Figure 9** then shows that the arrival intensity becomes low as $\eta$ increases, and it almost vanishes for a sufficiently large $\eta$ at which harvesting the fish is not beneficial compared with the utility gained by possibly sustaining the fish population; in **Figure 9(d)**, the fish will be harvested only near the terminal time if the fish population is small, with which sustaining the population will no longer be meaningful for the representative angler and he/she attempt to increase the value function only by the harvesting.

For **Figures 10 and 11**, we observe a monotone relationship that increasing $\psi$ implies more overestimating the terminal biomass, and hence results in a larger value function $\Phi$. For a positive $\psi$ corresponding to convex $\rho$, increasing $\psi$ leads to equilibrium controls with more intensively arriving at the river owing to the overestimation of the terminal biomass. For a negative $\psi$ corresponding to concave $\rho$, the harvesting is less active. We also conducted numerical experiments for different years using growth models for each year (**Figure A4**), suggesting that the shape of the average growth qualitatively affects the computed equilibrium control, for instance, the case in 2024 is qualitatively different from the other cases owing to the concave average growth while the other present sigmoidal growth curves.

Finally, we investigate the nonlinear expectation $G = G(t,x)$ in (27), which, in the proposed model, is the predicted terminal biomass conditioned on the current state $(t,X_t) = (t,x)$ for different values of $\eta$ and $\psi$. This $G$ is computable at each time step once the computed $g$ becomes available (see (27) which connects $G$ and $g$). **Figures 12 and 13** show the computed $G$ values for different $\eta$ and $\psi$ values, respectively. **Figure 12** suggests that $G$ smoothly increases as $\eta$ increases, which is due to more weighting of the sustainability concern than of the harvesting benefit, resulting in a smaller arrival intensity. **Figure 13** suggests that the transition of the computed $G$ for $\psi$ is again monotone such that a larger value of $\psi$ results in a more optimistic prediction of the sustainability of the fish in view of the terminal biomass.

The computed results for the value function, equilibrium control, and predicted biomass obtained in this study would be useful in applications because they provide information on how fish should be harvested and their consequences simultaneously, even for cases where dynamic programming does not apply. From the standpoint of the management of *Ayu*, the computational results suggest that the harvesting behavior of the representative angler is affected by the convex/concave nature of the terminal utility and growth curve model. This implies that efforts should be made for the data collection to accurately evaluate the spawning size (the nonlinear expectation in the terminal utility). Finding a suitable growth curve model, as we addressed in this paper, is therefore critical because it contributes to both the harvesting profit and terminal utility.

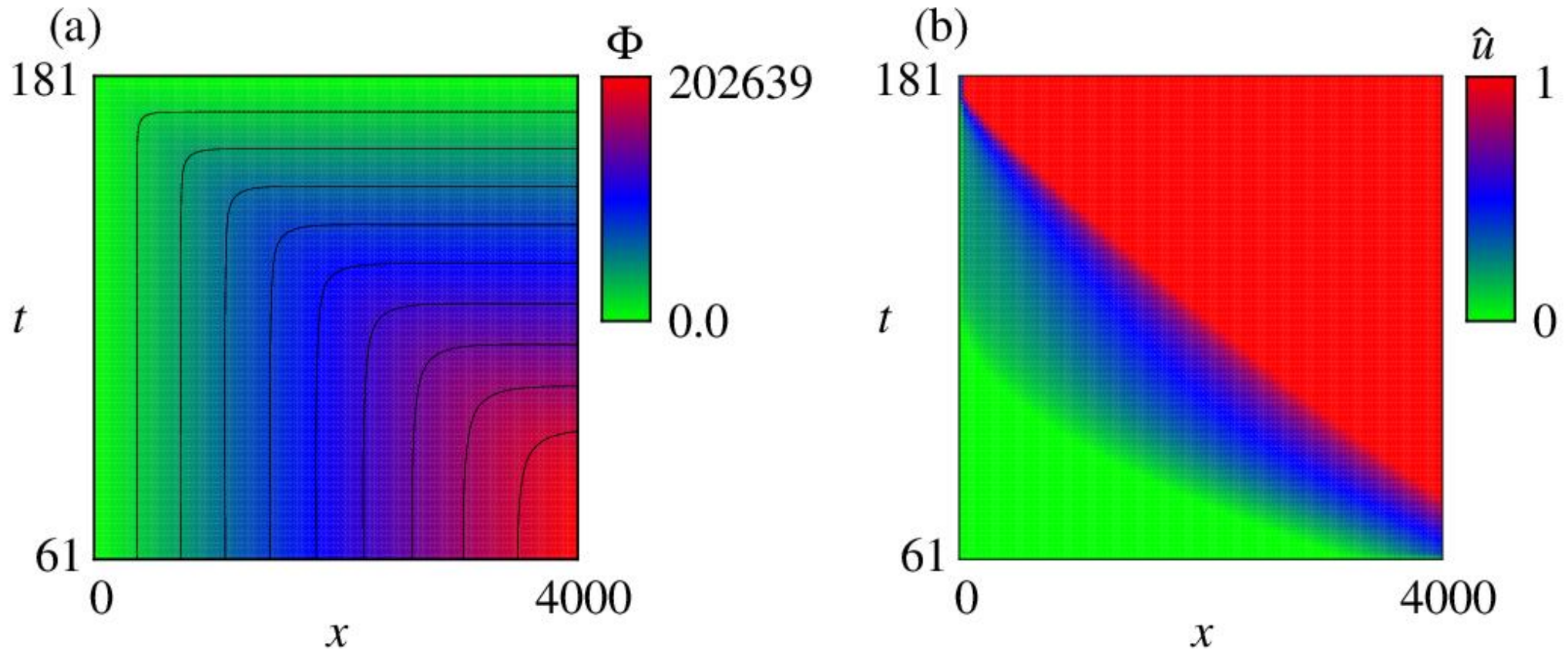

**Figure 6.** Computational results of (a) value function $\Phi$ and (b) the corresponding equilibrium control $\hat{u}$ when no terminal utility $\eta = 0$ is present (this case is independent from $\psi$ ).

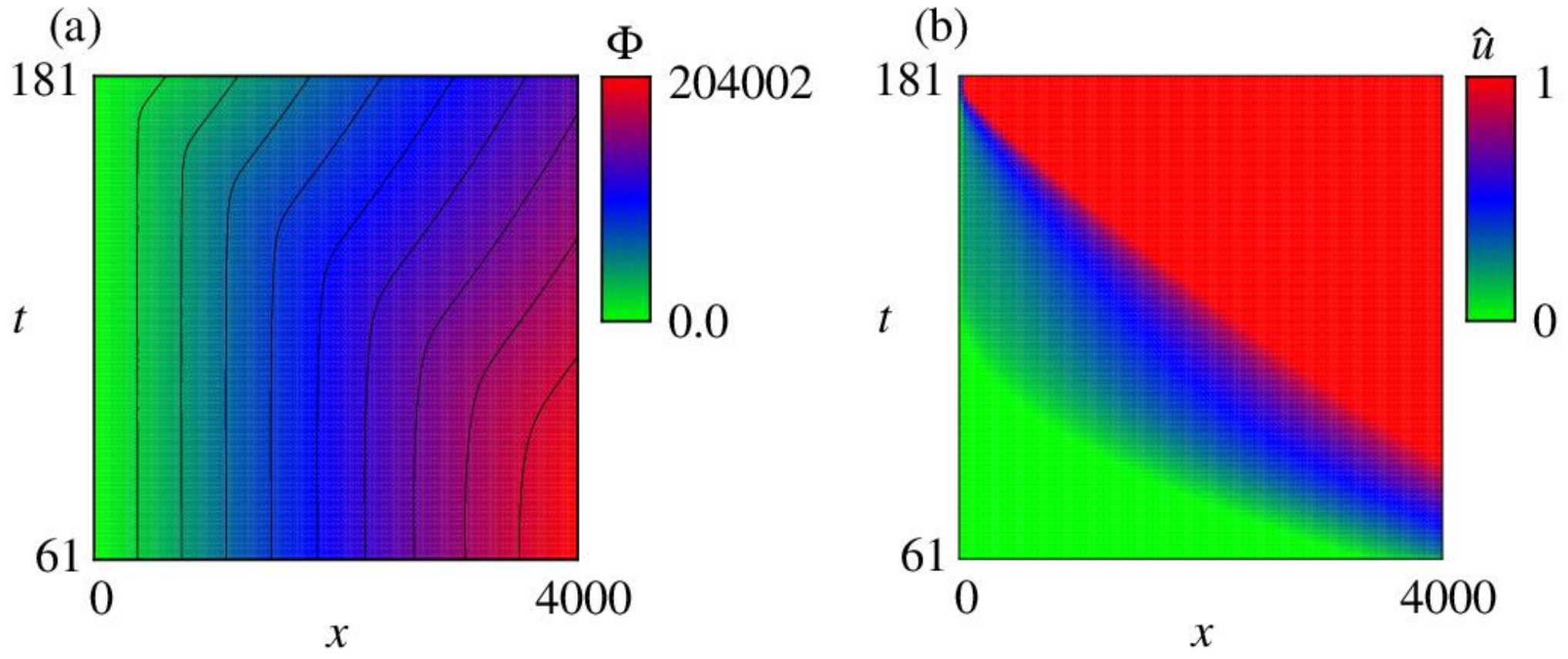

**Figure 7.** Computational results of (a) value function $\Phi$ and (b) the corresponding equilibrium control $\hat{u}$ for the benchmark case ( $\eta = 0.6$ and $\psi = 0$ ).

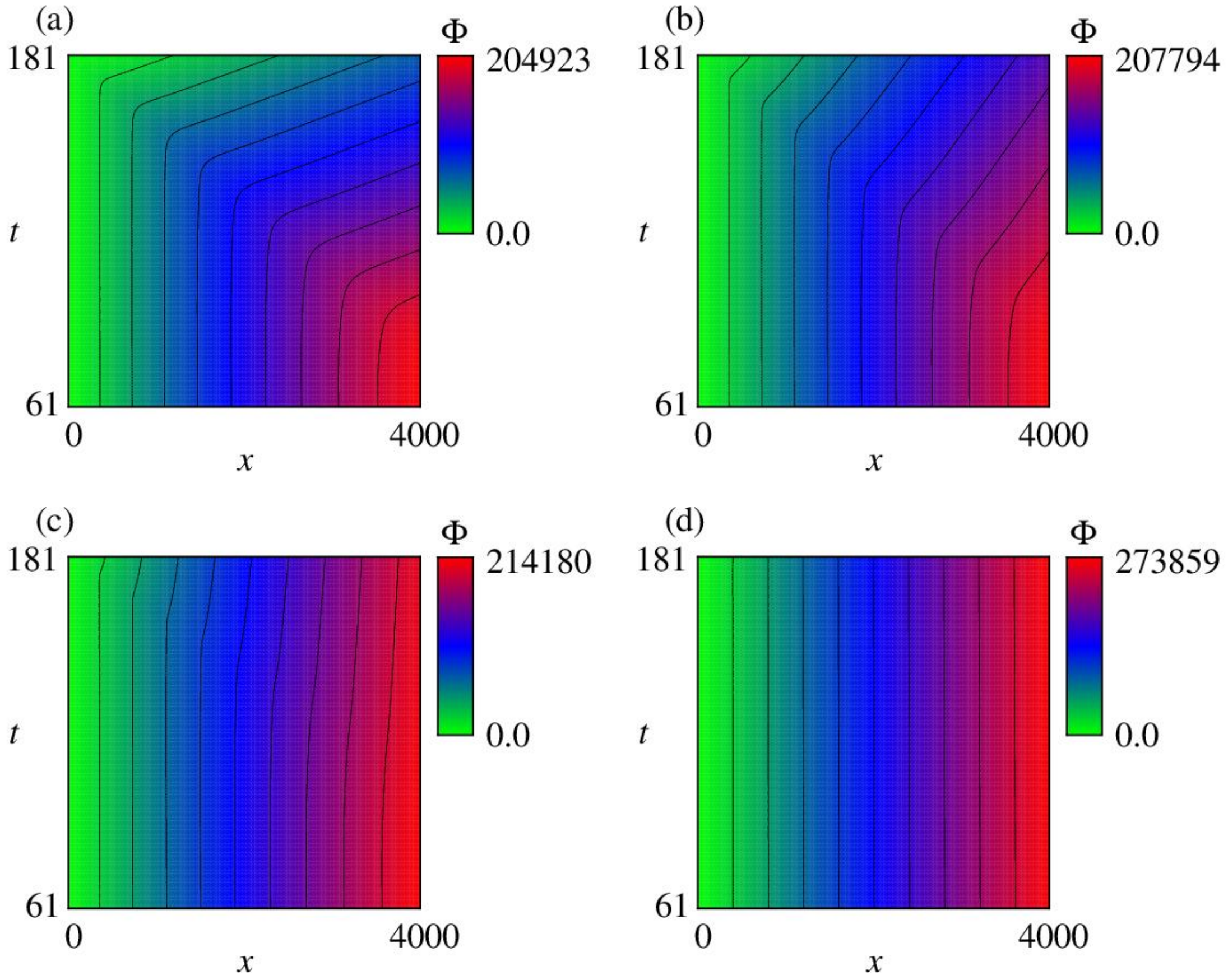

**Figure 8.** Computed value functions $\Phi$ with $\psi = 1.5$ for different values of $\eta$ : (a) $\eta = 0.3$ , (b) $\eta = 0.6$ , (c) $\eta = 0.9$ , and (d) $\eta = 1.2$ .

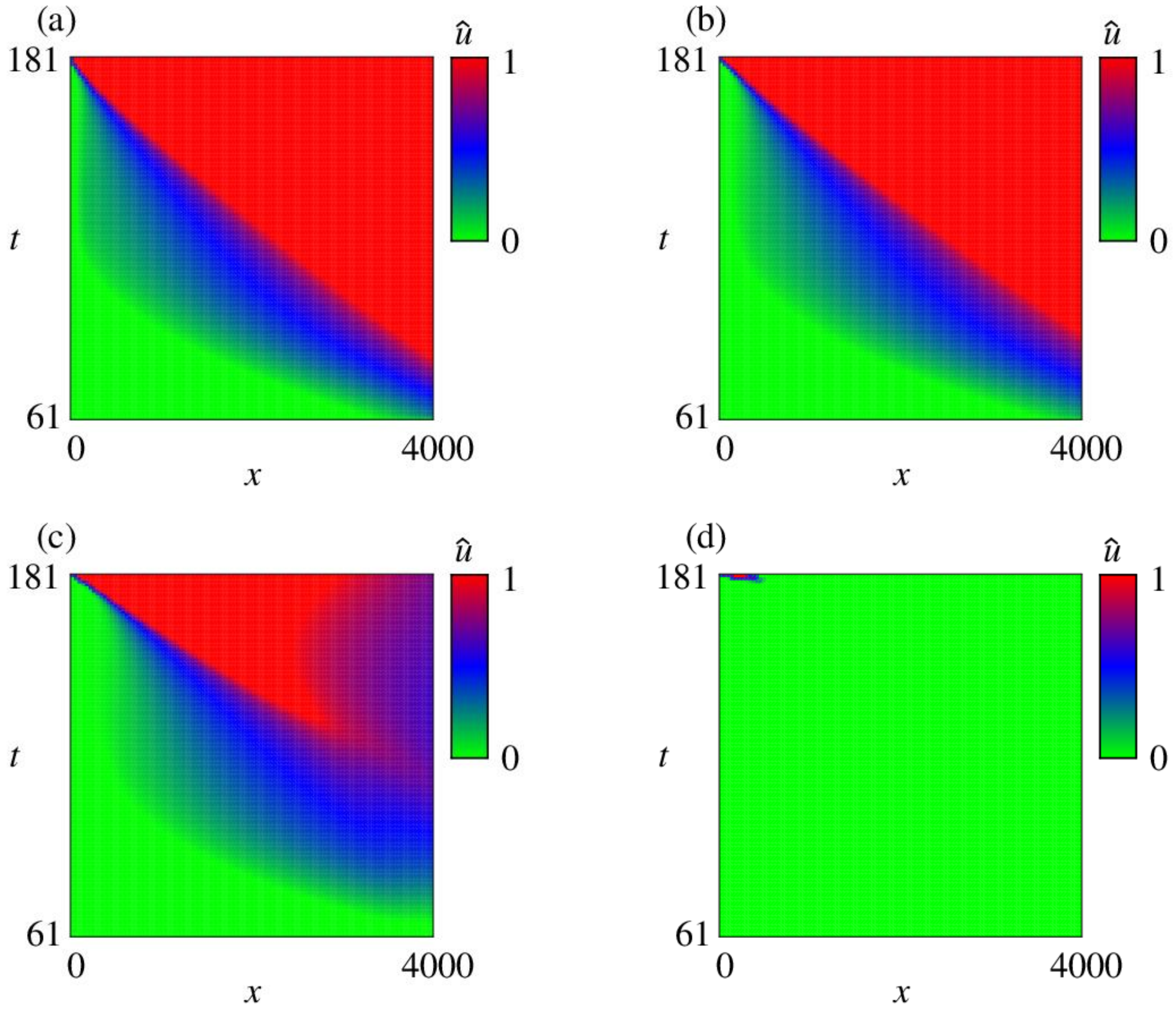

**Figure 9.** Computed equilibrium controls $\hat{u}$ with $\psi = 1.5$ for different values of $\eta$: (a) $\eta = 0.3$, (b) $\eta = 0.6$, (c) $\eta = 0.9$, and (d) $\eta = 1.2$.

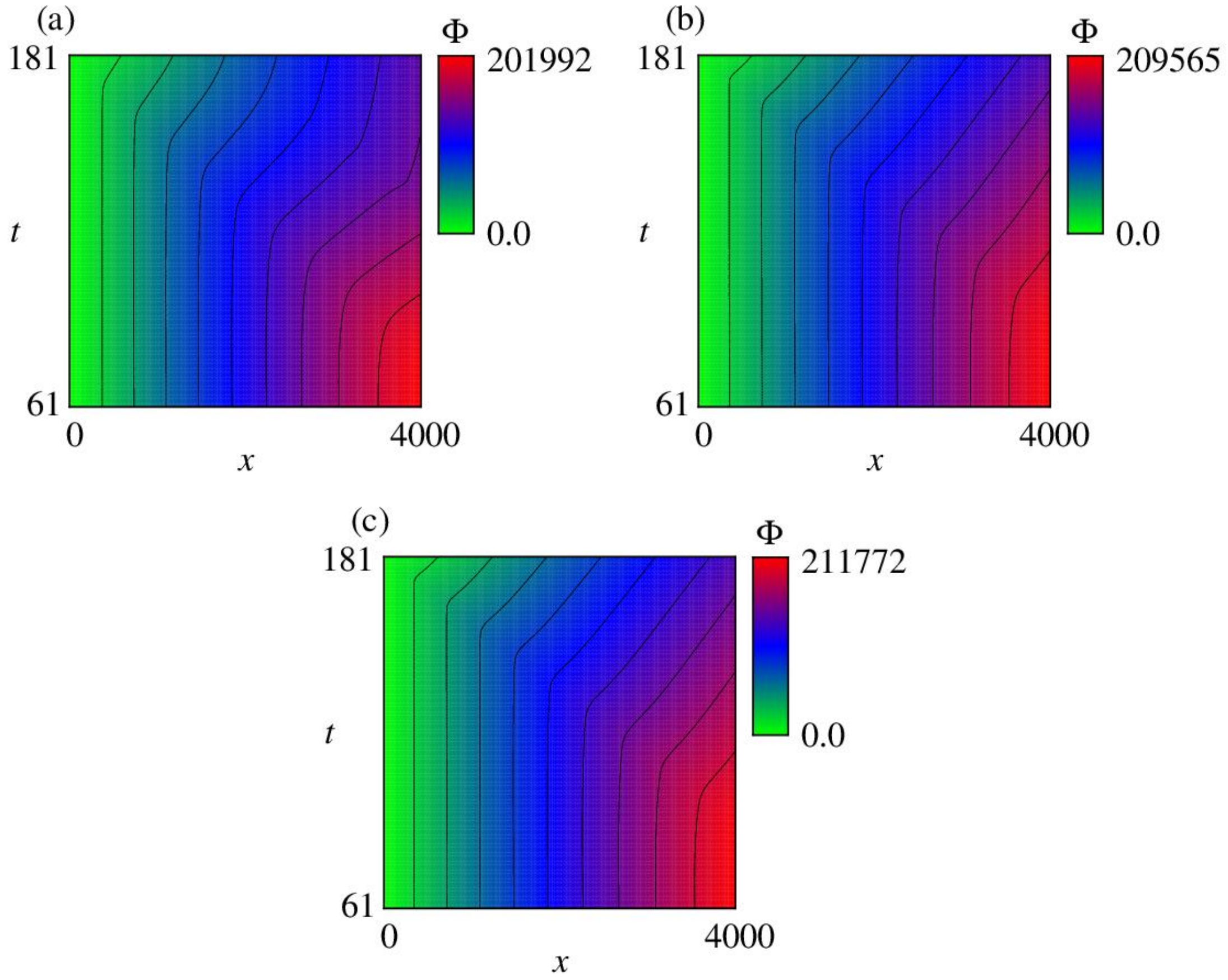

**Figure 10.** Computed value functions $\Phi$ with $\eta = 0.6$ for different values of $\psi$ : (a) $\psi = -0.75$, (b) $\psi = 2.5$, and (c) $\psi = 4.0$.

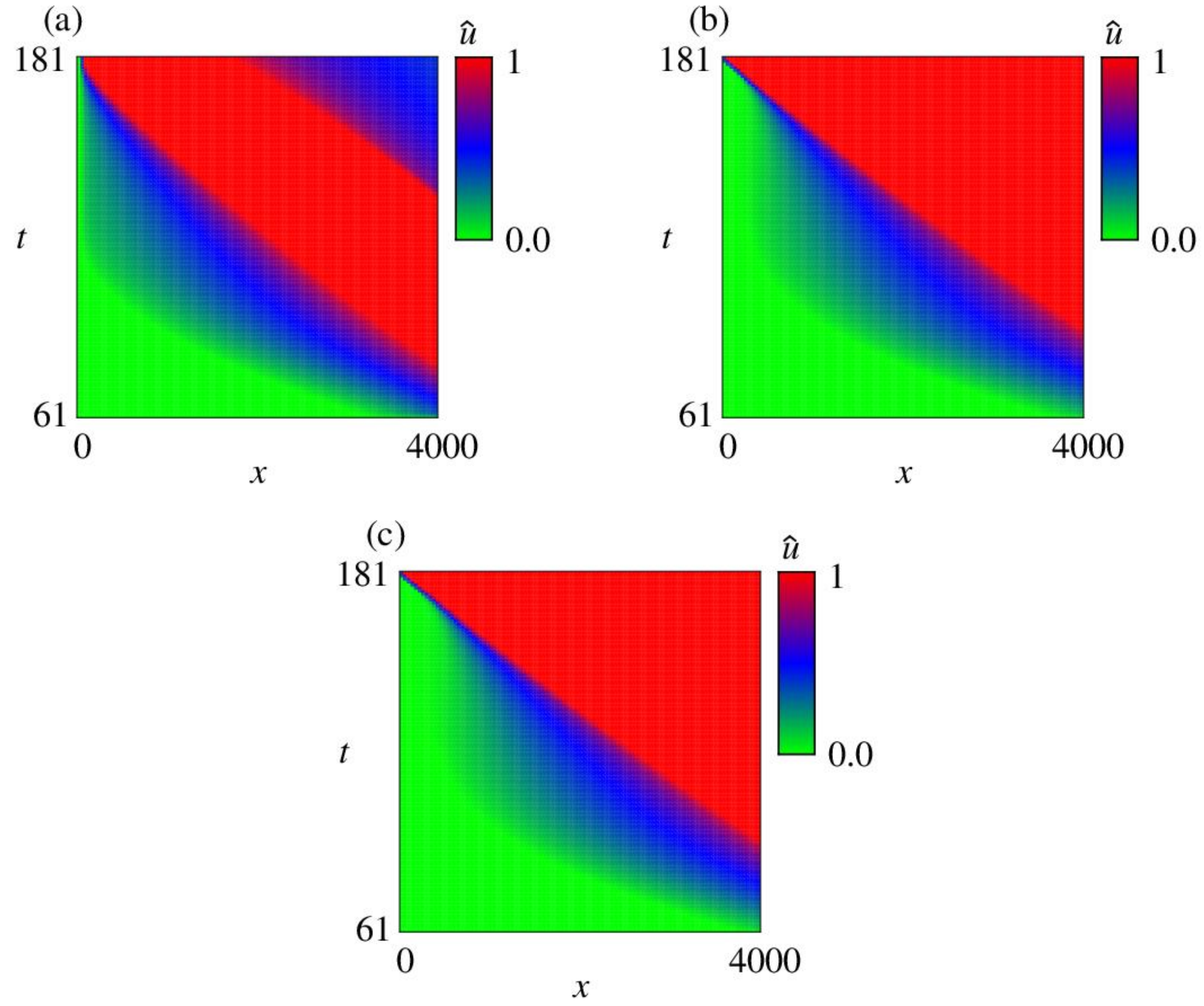

**Figure 11.** Computed equilibrium controls $\hat{u}$ with $\eta = 0.6$ for different values of $\psi$ : (a) $\psi = -0.75$, (b) $\psi = 2.5$, and (c) $\psi = 4.0$.

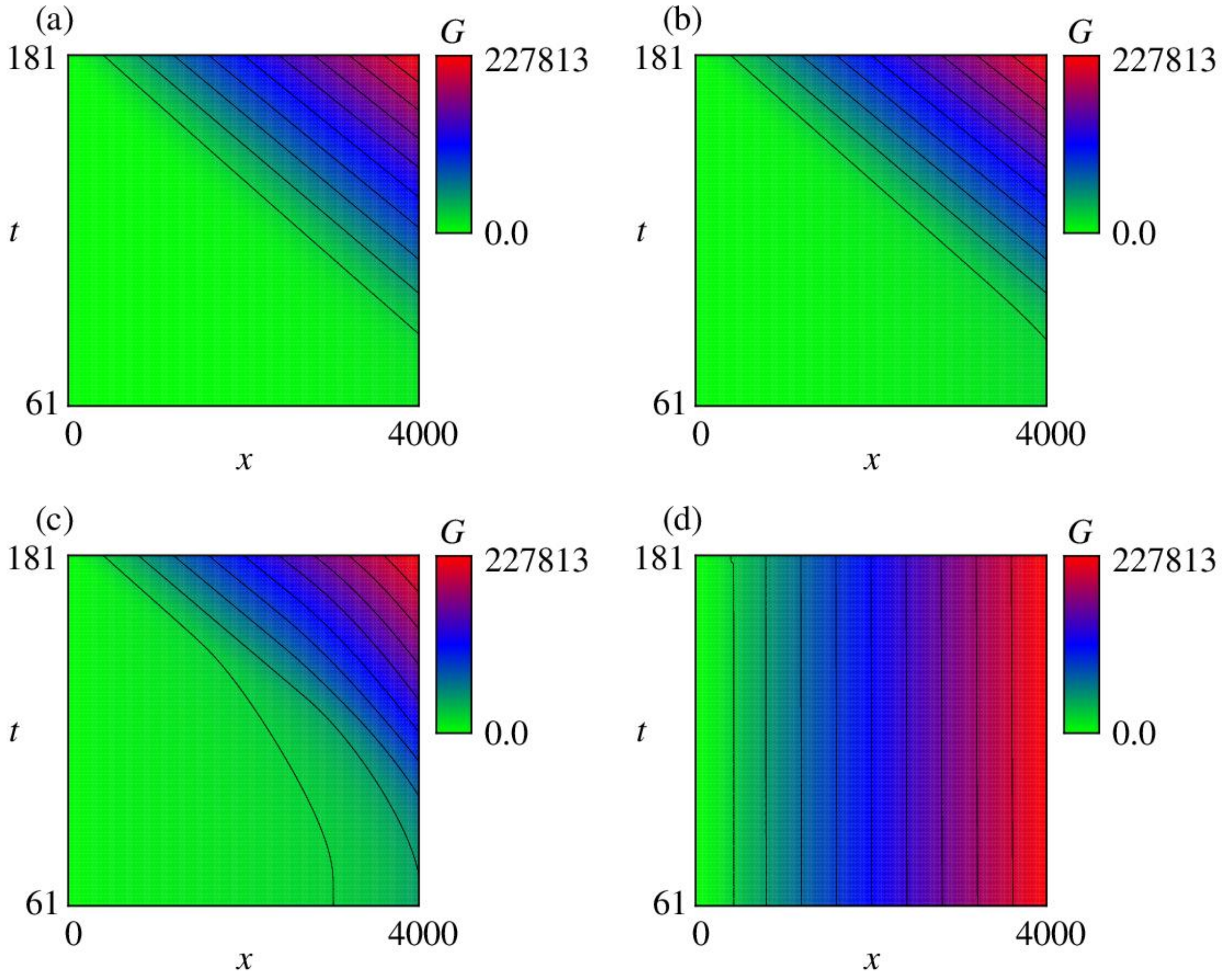

**Figure 12.** Computed nonlinear expectations $G$ with $\psi = 1.5$ for different values of $\eta$: (a) $\eta = 0.3$, (b) $\eta = 0.6$, (c) $\eta = 0.9$, and (d) $\eta = 1.2$.

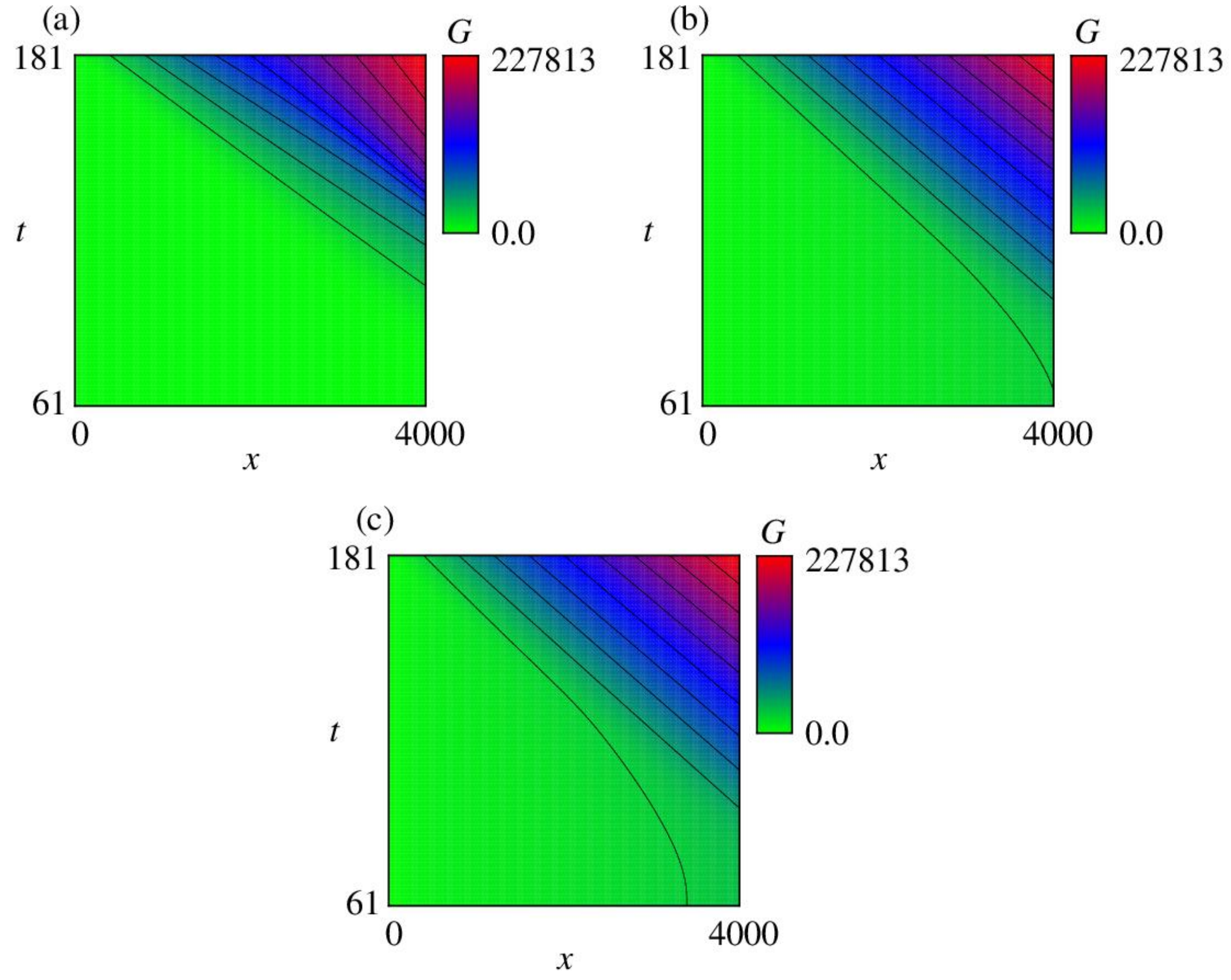

**Figure 13.** Computed nonlinear expectation $G$ with $\eta = 0.6$ for different values of $\psi$ : (a) $\psi = -0.75$, (b) $\psi = 2.5$, and (c) $\psi = 4.0$.

## 5. Conclusion

We applied a growth curve model with a size spectrum to the latest data set of *Ayu* in the Hii River system in Japan. The results demonstrated that the size spectrum was not negligible for the target fish. We found that logistic models fit the data well and that collecting data during the early growth stage (May to June) is crucial for investigating the growth of *Ayu* at the study site. We also developed a time-inconsistent stochastic control theory for fisheries management to account for the size spectrum. The computational results suggested that our control-theoretic approach simultaneously provide the value function, equilibrium control, and terminal nonlinear expectation, providing effective information toward sustainable fisheries management.

Fisheries management under multispecies interactions is an unresolved issue because other inland fishery resources in addition to *Ayu* are present in a river, which can be addressed theoretically by considering food web dynamics in a control problem; however, in this case, one will need to handle a high-dimensional nonlinear PDE system. This will be a theoretically and computationally more complicated issue than that studied in this paper but would be an important challenge toward the development of a more realistic fisheries management model. Another unresolved issue is the development of weak solution theory for time-inconsistent control problems, where the optimality equation, such as the HJB system, does not necessarily admit smooth solutions. Stronger nonlinearity than that arising in the classical control theory, such as that found in our extended HJB system, hinders the use of a well-established theory of viscosity solutions to the HJB equations [83]. This is also important for studying the convergence of numerical solutions to suitable solutions in extended HJB systems. The development of long-term [84] and mean-field models [85] specialized in inland fishery management is also an important future topic. While a model that adaptively incorporates both the catch season and the closed season would be interesting, we leave this for future work. Finally, the authors are currently developing a model that accounts for the causal relationship between water and air temperatures and fish growth.

**Appendix of "Mathematical model for sustainable fisheries resource management accounting for size spectrum" by Yoshioka et al.**

## A. Auxiliary results

This section presents the auxiliary results of the logistic model with time-dependent growth rate: **Tables A2 through A4** show the monthly air temperature statistics for Yokota shown in **Figure 1**. **Figure A1** compares empirical and theoretical average body weight each year for the logistic model with constant growth rate $r$, where we also plot the corresponding average ± standard deviation curves. **Figure A2** compares empirical and theoretical average body weight each year for the logistic model with time-dependent growth rate $r_t = r_0 + r_1 t$, where we also plot the corresponding average ± standard deviation curves. The fitted parameter values are listed in **Table A1**. **Figure A3** shows the river water temperatures (HOBO U-20) measured at Yumura, Kisuki, and Shinmitoya in **Figure 1**.

**Figure A4** shows the computed equilibrium controls $\hat{u}$ with $\eta = 0.6$ and $\psi = -0.75$ for different years; the growth model is a logistic model with a time-dependent growth rate. The profile of $\hat{u}$ depends on the growth model and has qualitatively similar shapes, except for that in 2024, where the growth is somewhat different from the other years, as shown in **Figure A2**.

**Table A1.** Fitted parameter values and the minimized Err (MinErr): logistic model with time-dependent growth rate.

| | 2017 | 2018 | 2019 | 2023 | 2024 | 2025 |
|---|---|---|---|---|---|---|
| $\alpha$ (-) | 8.47.E+00 | 9.59.E+00 | 9.60.E+00 | 6.18.E+00 | 1.01.E+01 | 8.36.E+00 |
| $\beta$ (g) | 1.29.E+01 | 1.01.E+01 | 1.17.E+01 | 1.26.E+01 | 1.52.E+05 | 6.83.E+00 |
| $f_0$ (g) | 1.76.E-01 | 1.16.E-01 | 1.17.E-01 | 3.60.E-01 | 3.24.E-05 | 1.99.E-01 |
| $W_0$ (g) | 1.92.E+01 | 1.12.E+01 | 1.31.E+01 | 2.81.E+01 | 4.99.E+01 | 1.14.E+01 |
| $r_0$ (1/day) | 1.63.E-02 | 2.50.E-02 | 2.14.E-02 | 1.42.E-02 | 0.00.E+00 | 2.70.E-02 |
| $r_1$ (1/day$^2$) | 4.61.E-05 | 4.72.E-06 | 2.24.E-05 | 2.34.E-04 | 1.75.E-05 | 6.39.E-04 |
| MinErr (g$^2$) | 3.98.E+01 | 5.89.E+01 | 4.31.E+01 | 1.61.E+02 | 8.03.E+01 | 8.57.E+01 |

**Table A2.** Mean air temperature at Yokota. Red: 1st place in the year and Magenta: second place in the year.

| | 2017 | 2018 | 2019 | 2020 | 2021 | 2022 | 2023 | 2024 | 2025 |
|---|---|---|---|---|---|---|---|---|---|
| March | 3.7 | 6.7 | 5.6 | 6.5 | 7.6 | 6.7 | 8.3 | 5.4 | 5.7 |
| April | 11.8 | 12.3 | 9.7 | 8.2 | 10.8 | 12.0 | 11.8 | 14.0 | 11.5 |
| May | 16.8 | 16.5 | 16.6 | 16.6 | 16.2 | 15.8 | 16.3 | 16.2 | 16.2 |
| June | 18.4 | 19.9 | 19.5 | 21.3 | 20.2 | 21.2 | 20.7 | 21.0 | 22.3 |
| July | 25.1 | 25.8 | 23.3 | 22.4 | 24.6 | 24.7 | 25.4 | 26.1 | 26.8 |
| August | 24.7 | 25.5 | 25.3 | 26.2 | 24.1 | 25.5 | 26.7 | 26.4 | 26.0 |
| September | 18.8 | 19.8 | 21.7 | 20.5 | 20.7 | 21.6 | 22.7 | 24.0 | 23.1 |
| October | 14.4 | 13.8 | 15.4 | 13.1 | 14.9 | 13.8 | 13.5 | 16.7 | 16.3 |

**Table A3.** Same with **Table A2** but the data present the maximum in daily air temperature.

| | 2017 | 2018 | 2019 | 2020 | 2021 | 2022 | 2023 | 2024 | 2025 |
|---|---|---|---|---|---|---|---|---|---|
| March | 9.8 | 14.3 | 11.7 | 12.8 | 14.8 | 13.5 | 15.9 | 10.6 | 11.9 |
| April | 18.2 | 19.9 | 16.5 | 14.9 | 18.5 | 19.9 | 18.6 | 21.4 | 18.8 |
| May | 24.1 | 22.5 | 24.4 | 22.9 | 22.3 | 23.4 | 23.2 | 23.2 | 23.0 |
| June | 24.8 | 25.4 | 25.0 | 27.3 | 26.1 | 27.2 | 26.5 | 26.9 | 28.8 |
| July | 30.5 | 31.5 | 28.2 | 26.3 | 30.4 | 30.1 | 31.0 | 31.1 | 33.8 |
| August | 30.0 | 31.7 | 30.8 | 32.7 | 29.3 | 31.1 | 32.2 | 33.0 | 32.5 |
| September | 24.4 | 23.9 | 27.3 | 25.8 | 25.5 | 26.9 | 28.0 | 30.5 | 28.2 |
| October | 18.9 | 19.2 | 20.6 | 19.3 | 21.3 | 20.2 | 20.4 | 22.2 | 21.7 |

**Table A4.** Same with **Table A2** but the data present the minimum in daily air temperature.

| | 2017 | 2018 | 2019 | 2020 | 2021 | 2022 | 2023 | 2024 | 2025 |
|---|---|---|---|---|---|---|---|---|---|
| March | -1.3 | 0.0 | 0.3 | 0.9 | 0.9 | 1.0 | 1.7 | 0.7 | 0.3 |
| April | 4.9 | 5.2 | 3.4 | 1.5 | 3.2 | 5.0 | 5.1 | 7.5 | 3.9 |
| May | 10.0 | 10.5 | 8.8 | 10.6 | 10.6 | 8.7 | 10.5 | 9.8 | 9.9 |
| June | 12.2 | 14.9 | 14.6 | 15.9 | 15.1 | 15.9 | 16.0 | 16.1 | 16.9 |
| July | 21.1 | 20.9 | 19.7 | 19.6 | 19.9 | 21.0 | 20.8 | 22.0 | 21.3 |
| August | 20.5 | 20.5 | 21.0 | 21.0 | 20.3 | 21.3 | 22.8 | 21.6 | 21.3 |
| September | 14.1 | 16.3 | 17.5 | 16.3 | 16.6 | 17.3 | 19.1 | 19.5 | 19.2 |
| October | 10.4 | 9.0 | 10.7 | 8.1 | 10.0 | 8.4 | 8.2 | 12.1 | 12.1 |

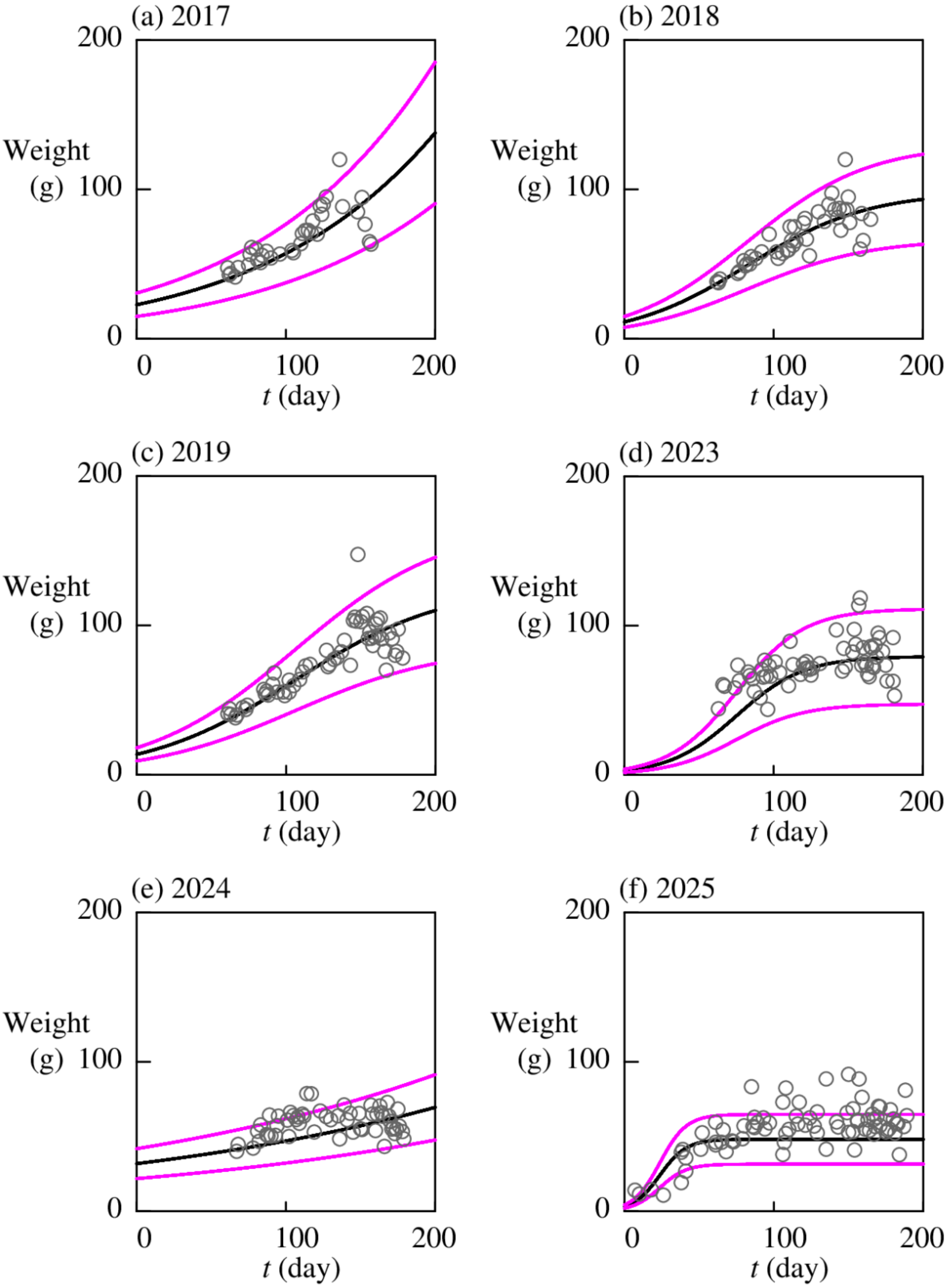

**Figure A1.** Comparison of average body weight each year: empirical (circles) and theoretical results of the logistic model (curves): average (black curves), average ± standard deviation (magenta curves).

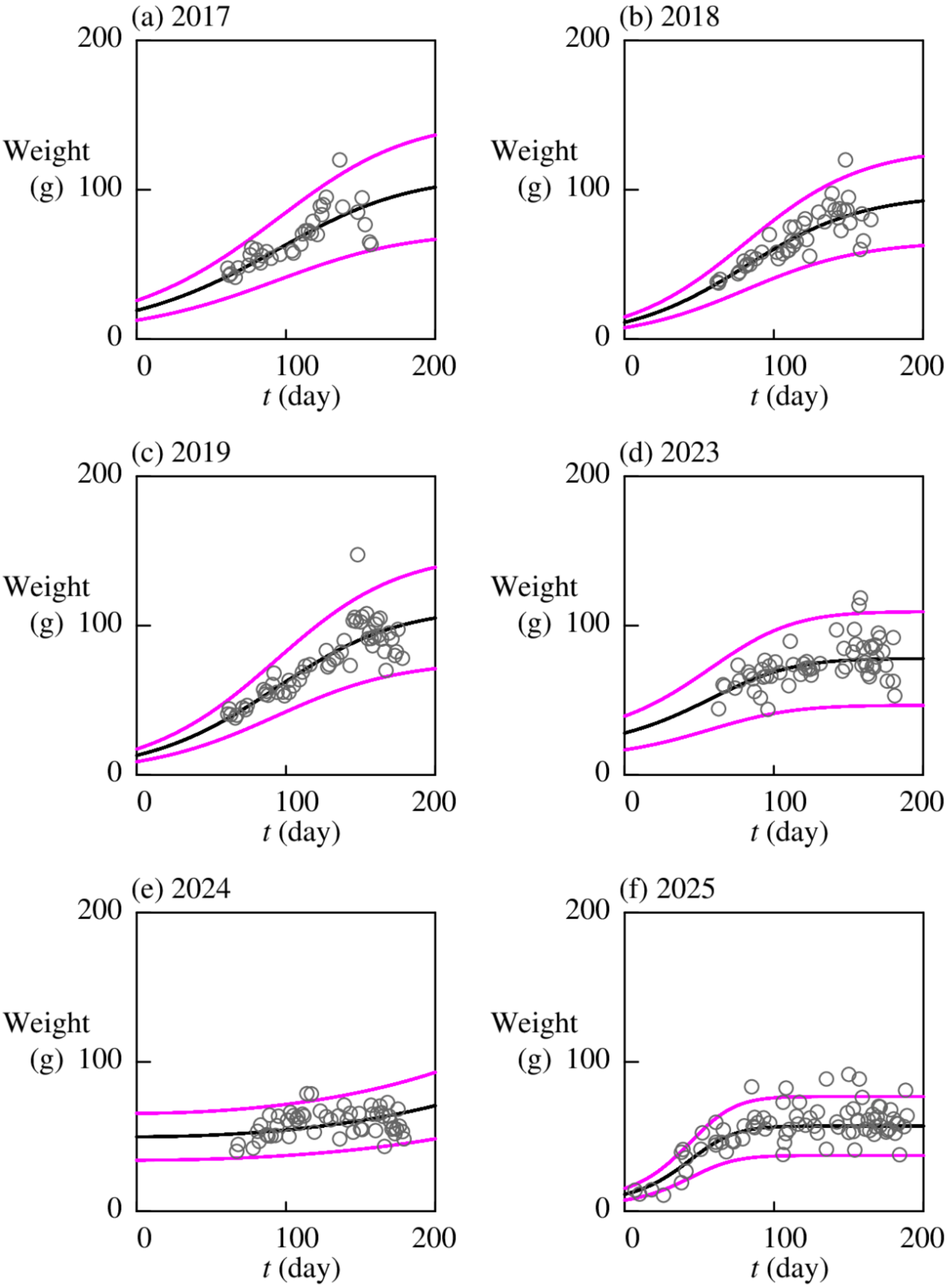

**Figure A2.** Comparison of average body weight each year: empirical (circles) and theoretical results of the logistic model with the time-dependent growth rate: average (black curves), average ± standard deviation (magenta curves).

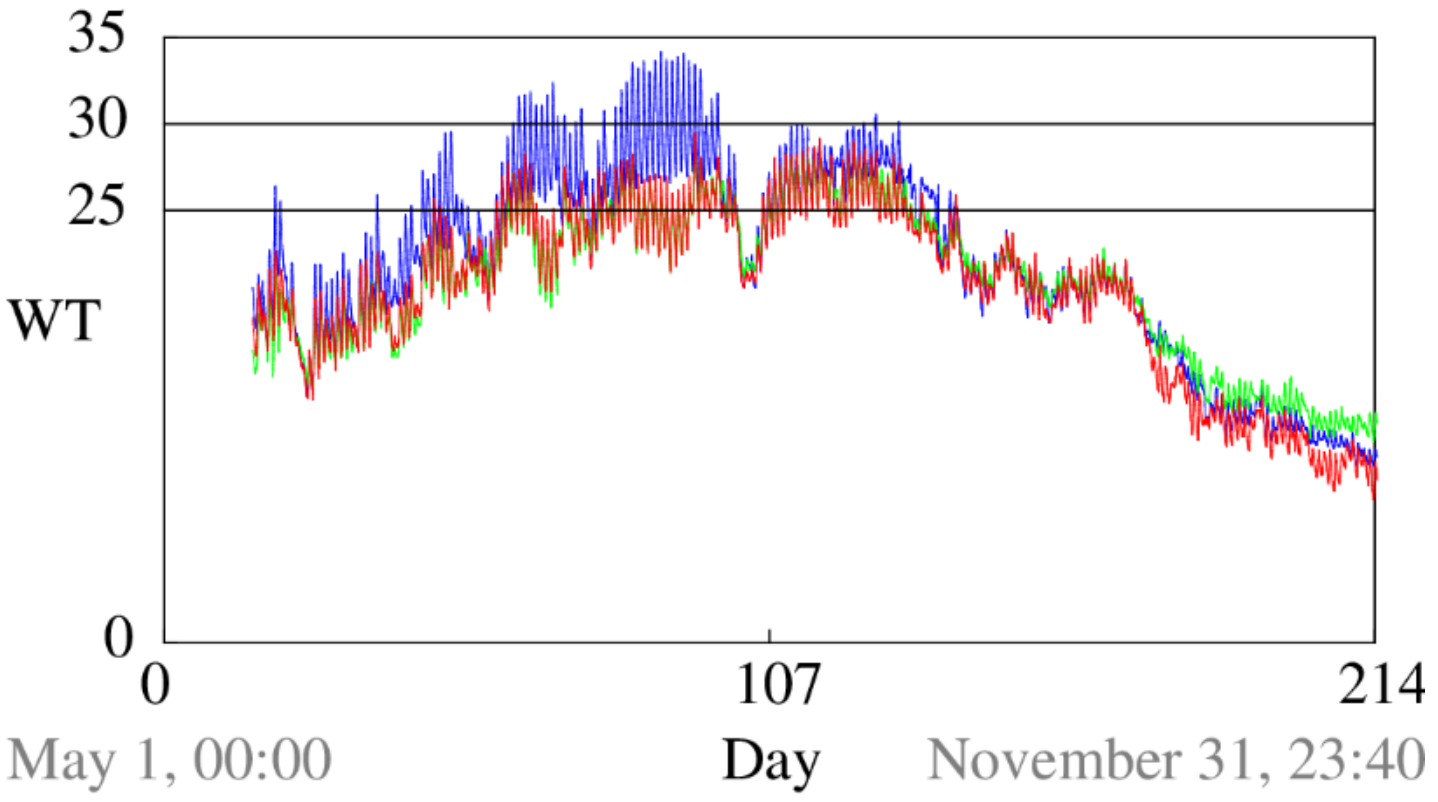

**Figure A3.** River water temperature (WT, °C) measured at Yumura (red), Kisuki (green), and Shinmitoya (blue) in **Figure 1**.

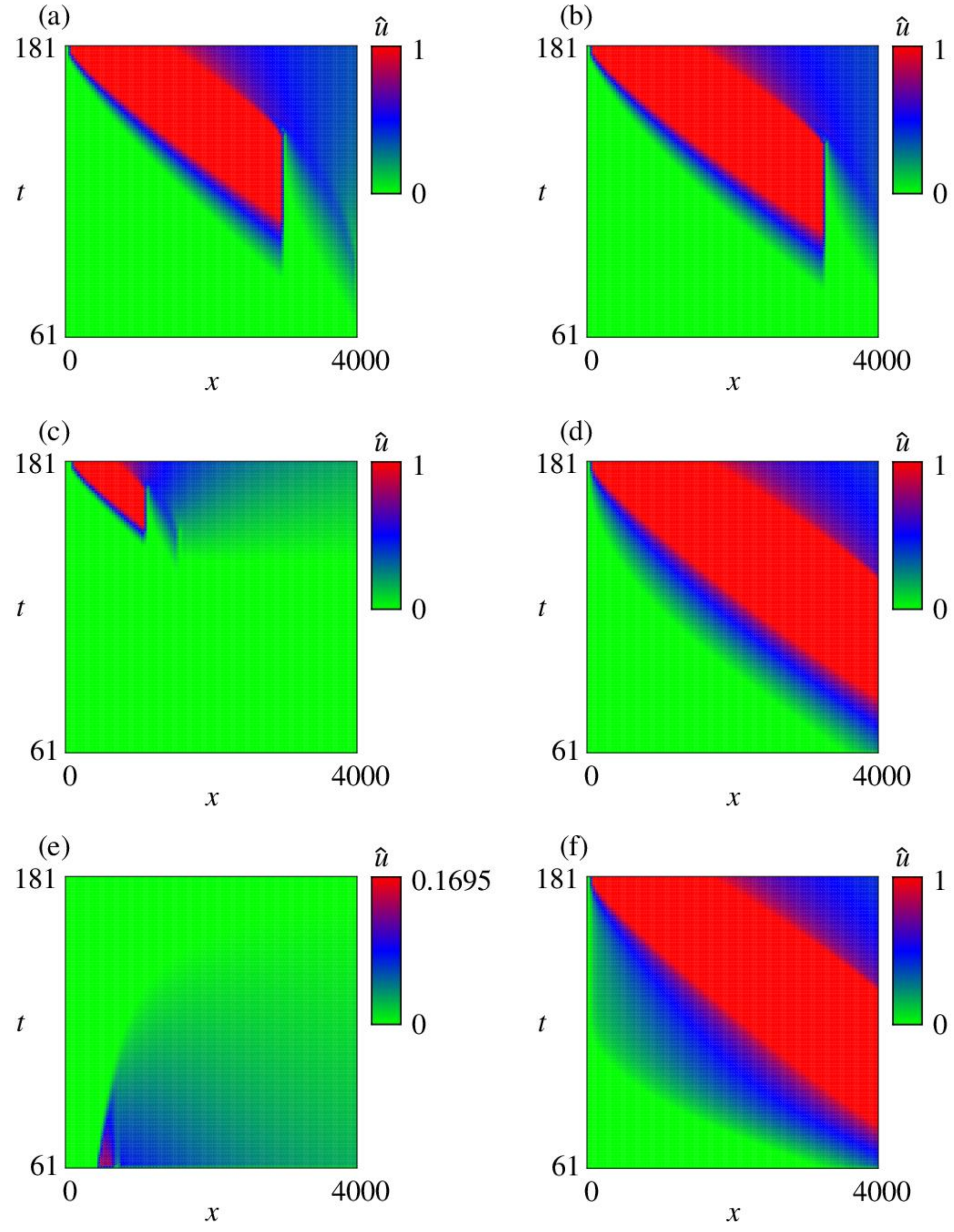

**Figure A4.** Computed equilibrium controls $\hat{u}$ with $\eta = 0.6$ and $\psi = -0.75$ for different years: (a) 2017, (b) 2018, (c) 2019, (d) 2023, (e) 2024, and (f) 2025.

## B. Numerical discretization

### B.1 Formulation

We present a finite difference method to discretize our extended HJB system because it does not seem to be analytically solvable. Our discretization strategy is simple. First, the state space is replaced with a grid, and then the extended HJB system is explicitly discretized on this grid.

Fix resolution parameters $N_t, N_x \in \mathbb{N}$. The temporal domain $[0,T]$ is uniformly discretized into a grid with vertices $s_i = i\Delta t$ ($i = 0,1,2,...,N_t$), where $\Delta t = T / N_t$. The spatial domain $[0,+\infty)$ is truncated to the bounded domain $[0,\bar{X}]$ with some $\bar{X} > 0$ and is discretized into a grid $x_j = j\Delta x$ ($j = 0,1,2,...,N_x$), where $\Delta x = T / N_x$. Vertices $P_{i,j}$ are defined as $(s_i, x_j)$ ($i = 0,1,2,...,N_t$, $j = 0,1,2,...,N_x$). The quantity $\Psi$ discretized at $P_{i,j}$ is expressed as $\Psi_{i,j}$, and if it does not depend on $i$ (resp., $j$), then it is simply written as $\Psi_j$ (resp., $\Psi_i$).

We explain a semi-discretized scheme with (28) being not discretized in the $w$ direction because its full discretization is completed by approximating the PDF $p_{W_t}$ (and hence $p_K$ owing to (6)) using suitable weighted Dirac deltas. We write $g(\cdot)$ discretized at $P_{i,j}$ as $g_{i,j}(\cdot)$. The terminal conditions (24) and (30) are discretized as $\Phi_{\cdot,j} = \eta \bar{W}_T x_j$ and $g_{\cdot,j}(\cdot) = \rho((\cdot)x_j)$, respectively. The boundary conditions (23) and (29) are implemented directly because they are constant. For simplicity of explanation, we assume the constant harvesting ($h(x) = \min\{\bar{h}, x\}$ with $\bar{h} > 0$) and possible complete extinction at catastrophe ($\kappa = 1$), and select $\bar{h} = \Delta x$ (see **Remark B1**). The other cases can be studied in an analogous way by linearly interpolating numerical solutions between each successive grid points.

With these preparations, we set the following recursion corresponding to (26) that is solved from $i = N_t$ to $i = 0$: at each $1 \le j \le N_x$,

$$\Phi_{i,j} = \Phi_{i+1,j} + \Delta t \max_{0\le\theta\le\bar{U}} \left\{ (\mathbb{L}_\theta \Phi)_{i,j} + \theta \bar{h} \bar{W}_{t_i} - \eta (\mathbb{L}_\theta G)_{i,j} + \eta \int_0^{+\infty} p_{W_T}(w) \lambda(g_{i+1,j}(w)) (\mathbb{L}_\theta g)_{i,j}(w) \mathrm{d}w \right\} \quad (38)$$

with

$$(\mathbb{L}_\theta \Phi)_{i,j} = \theta(\Phi_{i+1,j-1} - \Phi_{i+1,j}) + (d + k\theta^\gamma)(\Phi_{i+1,0} - \Phi_{i+1,j}), \quad (39)$$

$$\begin{aligned}(\mathbb{L}_\theta G)_{i,j} &= \theta \int_0^{+\infty} p_{W_T}(w) \left\{ \rho^{(-1)}(g_{i+1,j-1}(w)) - \rho^{(-1)}(g_{i+1,j}(w)) \right\} \mathrm{d}w \\ &\quad + (d + k\theta^\gamma) \int_0^{+\infty} p_{W_T}(w) \left\{ \rho^{(-1)}(g_{i+1,0}(w)) - \rho^{(-1)}(g_{i+1,j}(w)) \right\} \mathrm{d}w\end{aligned}, \quad (40)$$

and

$$(\mathbb{L}_\theta g)_{i,j}(w) = \theta(g_{i+1,j-1}(w) - g_{i+1,j}(w)) + (d + k\theta^\gamma)(g_{i+1,0}(w) - g_{i+1,j}(w)). \quad (41)$$

We set

$$\hat{\theta}_{i,j} = \arg\max_{0\le\theta\le\bar{U}} \left\{ (\mathbb{L}_\theta \Phi)_{i,j} + \theta \bar{h} \bar{W}_{t_i} - \eta (\mathbb{L}_\theta G)_{i,j} + \eta \int_0^{+\infty} p_{W_T}(w) \lambda(g_{i+1,j}(w)) (\mathbb{L}_\theta g)_{i,j}(w) \mathrm{d}w \right\}. \quad (42)$$

We can select such a $\hat{\theta}_{i,j}$ because the right-hand side of (42) is the maximizer of a continuous function in a closed interval. For each $1 \le i \le N_t - 1$, we set the discretization of (28) as follows: for each $1 \le j \le N_x$,

$$g_{i,j}(w) = g_{i+1,j}(w) + \Delta t \left( \mathbb{L}_{\hat{\theta}_{i,j}} g \right)_{i,j}(w). \tag{43}$$

In summary, we can compute (38) and (43) with (42), which can be performed explicitly in time because the right-hand sides of (38), (42), and (43) are available at time step $i$.

***Remark B1.*** We close this subsection with a remark on why the setting explained above suffices in **Section 4**. First, owing to the nonlocal structure in our extended HJB system such that the information at $(t,x)$ is formally constructed based on that at $(t, x-h(x))$; if $h(x) = \min\{x, \bar{h}\}$ and our interest is only for $x = \bar{h}i$ ($i = 0,1,2,...$), then selecting $\Delta x = \bar{h}$ suffices. Second, the choice of $\bar{X}$ essentially depends on the population size of each target problem. For each fixed grid, a numerical solution (i.e., $\Phi_{i,j}$ or $g_{i,j}(w)$) at $\mathrm{P}_{i,j}$ is formally based on those at $\mathrm{P}_{i+1,j-1}$ and $\mathrm{P}_{i+1,0}$, implying that numerical solutions can be computed from $i = N_t$ to $i = 0$ and from $j = 0$ to $j = N_x$; hence, domain-truncation error will not arise at $i = N_x$.

### B.2 Computational stability

We analyze the stability of the (semi-discretized) numerical solutions to our finite difference method for a special case. We have the following proposition, stating that numerical solutions are nonnegative and uniformly bounded if $\Delta t$ is suitably small. The upper bound of $\Delta t$ is independent of that of $\Delta x$. Condition (45) implies that $\rho$ does not increase rapidly for a large $x > 0$.

***Proposition B1***

*Assume that $\rho : [0,+\infty) \to [0,+\infty)$ is convex and strictly increasing. Assume further that $\Delta t$ is sufficiently small such that*

$$\Delta t \in \left( 0, \frac{1}{\bar{U} + d + k\bar{U}^{\gamma}} \right) \tag{44}$$

*and the following integrability condition*

$$\int_0^{+\infty} p_{W_T}(w) \lambda\left(\rho\left(w\bar{X}\right)\right)\left(\rho\left(w\bar{X}\right) - \rho(0)\right) \mathrm{d}w < +\infty. \tag{45}$$

*Then, at each $i = 0,1,2,...,N_t$ and $j = 0,1,2,...,N_x$, the following inequalities hold true:*

$$0 \le \Phi_{i,j} \le \eta \bar{W}_T \bar{h} \bar{X} + \bar{\Phi}\left(N_t - i\right)\Delta t \tag{46}$$

*and*

$$\rho(0) \le g_{i,j}(w) \le \rho\left(w\bar{X}\right) \text{ for all } w > 0. \tag{47}$$

*Here, $\bar{\Phi}$ is a positive constant independent of $\Delta t, \Delta x$.*

***Remark B2.*** The integrability condition in (45) is not restrictive. For example, this condition is satisfied for the power case (15) because

$$\int_0^{+\infty} p_{W_T}(w)\lambda\left(\rho\left(w\bar{X}\right)\right)\left(\rho\left(w\bar{X}\right)-\rho(0)\right)\mathrm{d}w = C\int_0^{+\infty} w p_{W_T}(w)\mathrm{d}w < +\infty \tag{48}$$

with a constant $C > 0$.

## C. Technical results

***Proof of Proposition 1***

The proof here basically follows that of Theorem 3.5 in Desmettre and Steffensen [39], but a significant modification is necessary because of the difference in system dynamics and objective functions between their problem and ours. In this proof, the state variable $X$ driven by $\hat{u}$ is expressed as $\hat{X}$.

**Step 1: Feynman–Kac representation for $g$**

Because $g(\cdot,\cdot,w)$ grows polynomially for all $w > 0$, applying Itô's formula (Theorem 1.14 in Øksendal and Sulem [43]) to $g$ along with its regularity assumption and (28) and (30) yields the Feynman–Kac representation

$$g(t,x,w) = \mathbb{E}^{t,x}\left[\rho\left(w\hat{X}_T\right)\right] \text{ for all } w > 0. \tag{49}$$

**Step 2: HJB-like equation**

The HJB-like equation (26) can be rewritten as follows: for any $0 \le t < T$ and $x > 0$,

$$\partial_t \Phi(t,x) + \left\{ \begin{array}{l} \mathbb{L}_\theta \Phi(t,x) + \theta h(x)\bar{W}_t - \eta\left(\partial_t G(t,x) + \mathbb{L}_\theta G(t,x)\right) \\ +\eta\int_0^{+\infty} p_{W_T}(w)\lambda\left(g(t,x,w)\right)\left(\partial_t g(t,x,w) + \mathbb{L}_\theta g(t,x,w)\right)\mathrm{d}w \end{array} \right\}_{\theta=\hat{\theta}(t,x)} = 0. \tag{50}$$

Substituting (28) into (50) yields the following equation without terms explicitly including $g$:

$$\partial_t \Phi(t,x) + \left\{\mathbb{L}_\theta \Phi(t,x) + \theta h(x)\bar{W}_t - \eta\left(\partial_t G(t,x) + \mathbb{L}_\theta G(t,x)\right)\right\}_{\theta=\hat{\theta}(t,x)} = 0. \tag{51}$$

Again, by the Itô's formula, for any $0 \le t \le T$ and $x > 0$, we have

$$\mathbb{E}^{t,x}\left[\Phi\left(T,\hat{X}_T\right)\right] = \Phi(t,x) + \mathbb{E}^{t,x}\left[\int_t^T \left(\partial_s + \mathbb{L}_{\hat{\theta}\left(s,\hat{X}_s\right)}\right)\Phi\left(s,\hat{X}_s\right)\mathrm{d}s\right] \tag{52}$$

and

$$\mathbb{E}^{t,x}\left[G\left(T,\hat{X}_T\right)\right] = G(t,x) + \mathbb{E}^{t,x}\left[\int_t^T \left(\partial_s + \mathbb{L}_{\hat{\theta}\left(s,\hat{X}_s\right)}\right)G\left(s,\hat{X}_s\right)\mathrm{d}s\right]. \tag{53}$$

We also have

$$
\begin{aligned}
\eta G(T,x) &= \eta\int_0^{+\infty} p_{W_T}(w)\rho^{(-1)}\left(g(T,x,w)\right)\mathrm{d}w \\
&= \eta\int_0^{+\infty} p_{W_T}(w)\rho^{(-1)}\left(\rho(wx)\right)\mathrm{d}w \\
&= \eta\left(\int_0^{+\infty} w p_{W_T}(w)\mathrm{d}w\right)x \\
&= \Phi(T,x)
\end{aligned}
. \tag{54}
$$

Substituting (51), (53), and (54) into (52) yields

$$
\begin{aligned}
&\mathbb{E}^{t,x}\left[\Phi\left(T,\hat{X}_T\right)\right] \\
&= \Phi(t,x) + \mathbb{E}^{t,x}\left[-\int_t^T \hat{\theta}\left(s,\hat{X}_s\right)h\left(\hat{X}_s\right)\bar{W}_s\mathrm{d}s + \eta\int_t^T\left(\partial_s + \mathbb{L}_{\hat{\theta}\left(s,\hat{X}_s\right)}\right)G\left(s,\hat{X}_s\right)\mathrm{d}s\right] \\
&= \Phi(t,x) - \mathbb{E}^{t,x}\left[\sum_{\substack{k=1,2,3,\ldots \\ t<\tau_k<T}} \bar{W}_{\tau_k} h\left(\hat{X}_{\tau_k -}\right)\right] + \eta\mathbb{E}^{t,x}\left[\int_t^T\left(\partial_s + \mathbb{L}_{\hat{\theta}\left(s,\hat{X}_s\right)}\right)G\left(s,\hat{X}_s\right)\mathrm{d}s\right], \\
&= \Phi(t,x) - \mathbb{E}^{t,x}\left[\sum_{\substack{k=1,2,3,\ldots \\ t<\tau_k<T}} \bar{W}_{\tau_k} h\left(\hat{X}_{\tau_k -}\right)\right] + \eta\left(\mathbb{E}^{t,x}\left[G\left(T,\hat{X}_T\right)\right] - G(t,x)\right) \\
&= \Phi(t,x) - \mathbb{E}^{t,x}\left[\sum_{\substack{k=1,2,3,\ldots \\ t<\tau_k<T}} \bar{W}_{\tau_k} h\left(\hat{X}_{\tau_k -}\right)\right] + \mathbb{E}^{t,x}\left[\Phi\left(T,\hat{X}_T\right)\right] - \eta G(t,x)
\end{aligned}
\tag{55}
$$

where

$$
\mathbb{E}^{t,x}\left[\int_t^T \hat{\theta}\left(s,\hat{X}_s\right)h\left(\hat{X}_s\right)\bar{W}_s\mathrm{d}s\right] = \mathbb{E}^{t,x}\left[\sum_{\substack{k=1,2,3,\ldots \\ t<\tau_k<T}} \bar{W}_{\tau_k} h\left(\hat{X}_{\tau_k -}\right)\right]. \tag{56}
$$

Rearranging (55) yields

$$
\Phi(t,x) = \mathbb{E}^{t,x}\left[\sum_{\substack{k=1,2,3,\ldots \\ t<\tau_k<T}} \bar{W}_{\tau_k} h\left(\hat{X}_{\tau_k -}\right)\right] + \eta G(t,x). \tag{57}
$$

Moreover, by (27) and (49), we have

$$
G(t,x) = \int_0^{+\infty} p_{W_T}(w)\rho^{(-1)}\left(g(t,x,w)\right)\mathrm{d}w = \int_0^{+\infty} p_{W_T}(w)\rho^{(-1)}\left(\mathbb{E}^{t,x}\left[\rho\left(w\hat{X}_T\right)\right]\right)\mathrm{d}w. \tag{58}
$$

Substituting (58) into (57) yields

$$
\Phi(t,x) = \mathbb{E}^{t,x}\left[\sum_{\substack{k=1,2,3,\ldots \\ t<\tau_k<T}} \bar{W}_{\tau_k} h\left(\hat{X}_{\tau_k -}\right)\right] + \eta\int_0^{+\infty} p_{W_T}(w)\rho^{(-1)}\left(\mathbb{E}^{t,x}\left[\rho\left(w\hat{X}_T\right)\right]\right)\mathrm{d}w = J(t,x,\hat{u}). \tag{59}
$$

In the rest of this proof, we prove (22) for any $u \in \mathbb{U}$ with which the proof is completed by (59).

**Step 3: Perturbed control and the equilibrium condition**

Choose a perturbed control $u_\delta \in \mathbb{U}$ as in (21) with $\hat{u}$ given by (33). The process $X$ driven by the control $u_\delta$ is denoted as $X^\delta$. As in "Step 3" of Desmettre and Steffensen [39], we obtain

$$J(t,x;u_\delta)=\mathbb{E}^{t,x}\left[J\left(t+\delta,X^\delta_{t+\delta};u_\delta\right)\right]+\mathbb{E}^{t,x}\left[\sum_{\substack{k=1,2,3,\ldots\\ t<\tau_k<t+\delta}}\bar{W}_{\tau_k}h\left(X^\delta_{\tau_k-}\right)\right]$$
$$-\eta\mathbb{E}^{t,x}\left[\int_0^{+\infty}p_{W_T}(w)\rho^{(-1)}\left(g^\delta\left(t+\delta,X^\delta_{t+\delta},w\right)\right)\mathrm{d}w\right] \quad , \tag{60}$$
$$+\eta\int_0^{+\infty}p_{W_T}(w)\rho^{(-1)}\left(\mathbb{E}^{t,x}\left[g^\delta\left(t+\delta,X^\delta_{t+\delta},w\right)\right]\right)\mathrm{d}w$$

where $g^\delta(t,x,w)=\mathbb{E}^{t,x}\left[\rho\left(wX^\delta_T\right)\right]$.

Because $u_\delta=u$ on $[t,t+\delta)$, if $X_t=X^\delta_t=x$, then

$$X_{t+\delta}=X^\delta_{t+\delta}. \tag{61}$$

Similarly, we have

$$\mathbb{E}^{t,x}\left[\sum_{\substack{k=1,2,3,\ldots\\ t<\tau_k<t+\delta}}\bar{W}_{\tau_k}h\left(X^\delta_{\tau_k-}\right)\right]=\mathbb{E}^{t,x}\left[\sum_{\substack{k=1,2,3,\ldots\\ t<\tau_k<t+\delta}}\bar{W}_{\tau_k}h\left(X_{\tau_k-}\right)\right]. \tag{62}$$

Moreover, because $u_\delta=\hat{u}$ on $[t+\delta,T]$,

$$J\left(t+\delta,X^\delta_{t+\delta};u_\delta\right)=J\left(t+\delta,X^\delta_{t+\delta};\hat{u}\right)=J\left(t+\delta,\hat{X}_{t+\delta};\hat{u}\right). \tag{63}$$

We also have

$$g^\delta\left(t+\delta,X^\delta_{t+\delta},\cdot\right)=g\left(t+\delta,X^\delta_{t+\delta},\cdot\right)=g\left(t+\delta,X_{t+\delta},\cdot\right). \tag{64}$$

Consequently, we can rewrite (60) as follows, where $u_\delta$ does not appear on the right side:

$$J(t,x;u_\delta)=\mathbb{E}^{t,x}\left[\Phi\left(t+\delta,X_{t+\delta}\right)\right]+\mathbb{E}^{t,x}\left[\sum_{\substack{k=1,2,3,\ldots\\ t<\tau_k<t+\delta}}\bar{W}_{\tau_k}h\left(X_{\tau_k-}\right)\right]$$
$$-\eta\begin{pmatrix}\mathbb{E}^{t,x}\left[\int_0^{+\infty}p_{W_T}(w)\rho^{(-1)}\left(g\left(t+\delta,X_{t+\delta},w\right)\right)\mathrm{d}w\right]\\ -\int_0^{+\infty}p_{W_T}(w)\rho^{(-1)}\left(\mathbb{E}^{t,x}\left[g\left(t+\delta,X_{t+\delta},w\right)\right]\right)\mathrm{d}w\end{pmatrix} \quad . \tag{65}$$
$$=\mathbb{E}^{t,x}\left[\Phi\left(t+\delta,X_{t+\delta}\right)\right]+\mathbb{E}^{t,x}\left[\sum_{\substack{k=1,2,3,\ldots\\ t<\tau_k<t+\delta}}\left(\int_0^{+\infty}wp_{W_{\tau_k}}(w)\mathrm{d}w\right)h\left(X_{\tau_k-}\right)\right]$$
$$-\eta\left(\mathbb{E}^{t,x}\left[G\left(t+\delta,X_{t+\delta}\right)\right]-\int_0^{+\infty}p_{W_T}(w)\rho^{(-1)}\left(\mathbb{E}^{t,x}\left[g\left(t+\delta,X_{t+\delta},w\right)\right]\right)\mathrm{d}w\right)$$

Here, we used (59).

By the HJB-like equation (26), we obtain the inequality

$$\partial_t\Phi(t,x)+\left\{\begin{array}{l}\mathbb{L}_\theta\Phi(t,x)+\theta h(x)\bar{W}_t-\eta\left(\partial_t G(t,x)+\mathbb{L}_\theta G(t,x)\right)\\ +\eta\int_0^{+\infty}p_{W_T}(w)\lambda\left(g(t,x,w)\right)\left(\partial_t g(t,x,w)+\mathbb{L}_\theta g(t,x,w)\right)\mathrm{d}w\end{array}\right\}_{\theta=u(t,x)}\leq 0 \tag{66}$$

since $\theta=u(t,x)$ is suboptimal for "max" in (26). By applying the Itô's formula to $\Phi,G,g(\cdot,\cdot,w)$ ($w>0$) in $(t,t+\delta)$ with $0<\delta<T-t$, we obtain

$$\mathbb{E}^{t,x}\left[\Phi\left(t+\delta,X_{t+\delta}\right)\right]-\Phi\left(t,x\right)=\delta\left(\partial_t\Phi\left(t,x\right)+\mathbb{L}_{u(t,x)}\Phi\left(t,x\right)\right)+o\left(\delta\right), \tag{67}$$

$$\mathbb{E}^{t,x}\left[G\left(t+\delta,X_{t+\delta}\right)\right]-G\left(t,x\right)=\delta\left(\partial_t G\left(t,x\right)+\mathbb{L}_{u(t,x)}G\left(t,x\right)\right)+o\left(\delta\right), \tag{68}$$

and

$$\mathbb{E}^{t,x}\left[g\left(t+\delta,X_{t+\delta},\cdot\right)\right]-g\left(t,x,\cdot\right)=\delta\left(\partial_t g\left(t,x,\cdot\right)+\mathbb{L}_{u(t,x)}g\left(t,x,\cdot\right)\right)+o\left(\delta\right). \tag{69}$$

Here, $o\left(\delta\right)$ represents terms such that $\delta^{-1}o\left(\delta\right)\to 0$ as $\delta\to +0$ . We also have

$$\begin{aligned}&\int_0^{+\infty}p_{W_T}\left(w\right)\rho^{(-1)}\left(\mathbb{E}^{t,x}\left[g\left(t+\delta,X_{t+\delta},w\right)\right]\right)\mathrm{d}w-\int_0^{+\infty}p_{W_T}\left(w\right)\rho^{(-1)}\left(g\left(t,x,w\right)\right)\mathrm{d}w\\&=\delta\int_0^{+\infty}p_{W_T}\left(w\right)\lambda\left(g\left(t,x,w\right)\right)\left(\partial_t g\left(t,x,w\right)+\mathbb{L}_\theta g\left(t,x,w\right)\right)\mathrm{d}w+o\left(\delta\right)\end{aligned}. \tag{70}$$

Now, (65) is rewritten as

$$\begin{aligned}&\mathbb{E}^{t,x}\left[\Phi\left(t+\delta,X_{t+\delta}\right)\right]\\&=J\left(t,x;u_\delta\right)-\mathbb{E}^{t,x}\left[\sum_{\substack{k=1,2,3,\ldots\\ t<\tau_k<t+\delta}}\bar{W}_{\tau_k}h\left(X_{\tau_k-}\right)\right]\\&+\eta\left(\mathbb{E}^{t,x}\left[G\left(t+\delta,X_{t+\delta}\right)\right]-\int_0^{+\infty}p_{W_T}\left(w\right)\rho^{(-1)}\left(\mathbb{E}^{t,x}\left[g\left(t+\delta,X_{t+\delta},w\right)\right]\right)\mathrm{d}w\right)\\&=J\left(t,x;u_\delta\right)-\delta u\left(t,x\right)h\left(x\right)\bar{W}_t\\&+\eta\left(\mathbb{E}^{t,x}\left[G\left(t+\delta,X_{t+\delta}\right)\right]-\int_0^{+\infty}p_{W_T}\left(w\right)\rho^{(-1)}\left(\mathbb{E}^{t,x}\left[g\left(t+\delta,X_{t+\delta},w\right)\right]\right)\mathrm{d}w\right)+o\left(\delta\right)\end{aligned}. \tag{71}$$

Combining (66)-(70) yields

$$\begin{aligned}&\mathbb{E}^{t,x}\left[\Phi\left(t+\delta,X_{t+\delta}\right)\right]-\Phi\left(t,x\right)+\delta u\left(t,x\right)h\left(x\right)\bar{W}_t\\&-\eta\left(\mathbb{E}^{t,x}\left[G\left(t+\delta,X_{t+\delta}\right)\right]-G\left(t,x\right)\right)\\&+\eta\left(\int_0^{+\infty}p_{W_T}\left(w\right)\rho^{(-1)}\left(\mathbb{E}^{t,x}\left[g\left(t+\delta,X_{t+\delta},w\right)\right]\right)\mathrm{d}w-G\left(t,x\right)\right)\mathrm{d}w\leq o\left(\delta\right)\end{aligned}, \tag{72}$$

which can be rearranged as

$$\begin{aligned}&\mathbb{E}^{t,x}\left[\Phi\left(t+\delta,X_{t+\delta}\right)\right]-\Phi\left(t,x\right)+\delta u\left(t,x\right)h\left(x\right)\bar{W}_t\\&-\eta\mathbb{E}^{t,x}\left[G\left(t+\delta,X_{t+\delta}\right)\right]+\eta\int_0^{+\infty}p_{W_T}\left(w\right)\rho^{(-1)}\left(\mathbb{E}^{t,x}\left[g\left(t+\delta,X_{t+\delta},w\right)\right]\right)\mathrm{d}w\leq o\left(\delta\right)\end{aligned}. \tag{73}$$

Substituting (71) into (73) yields

$$J\left(t,x,\hat{u}\right)=\Phi\left(t,x\right)\geq J\left(t,x,u_\delta\right)+o\left(\delta\right), \tag{74}$$

from which we obtain the desired inequality (22), and the proof is completed.

□

***Proof of Proposition B1***

We use an induction argument. Evidently, (46) and (47) are satisfied for $i=N_t$ and all $j=0,1,2,\ldots,N_x$ . Assume that (46) and (47) are satisfied for some $i+1\in\{1,2,3,\ldots,N_t\}$ and all $j=0,1,2,\ldots,N_x$ . First, from (41) and (43), for each $w>0$ , we obtain

$$
\begin{aligned}
g_{i,j}(w) &= g_{i+1,j}(w) + \Delta t\left(\mathbb{L}_{\hat{\theta}_{i,j}} g\right)_{i,j}(w) \\
&= g_{i+1,j}(w) + \Delta t\hat{\theta}_{i,j}\left(g_{i+1,j-1}(w) - g_{i+1,j}(w)\right) + \left(d + k\left(\hat{\theta}_{i,j}\right)^{\gamma}\right)\Delta t\left(g_{i+1,0}(w) - g_{i+1,j}(w)\right) \\
&= \left(1 - \Delta t\hat{\theta}_{i,j} - \left(d + k\left(\hat{\theta}_{i,j}\right)^{\gamma}\right)\Delta t\right) g_{i+1,j}(w) + \Delta t\hat{\theta}_{i,j} g_{i+1,j-1}(w) + \left(d + k\left(\hat{\theta}_{i,j}\right)^{\gamma}\right)\Delta t g_{i+1,0}(w)
\end{aligned} \quad . \tag{75}
$$

By (44) and the induction, for each $w > 0$, we have

$$
\begin{aligned}
g_{i,j}(w) &\geq \left(1 - \Delta t\hat{\theta}_{i,j} - \left(d + k\left(\hat{\theta}_{i,j}\right)^{\gamma}\right)\Delta t\right)\rho(0) + \Delta t\hat{\theta}_{i,j}\rho(0) + \left(d + k\left(\hat{\theta}_{i,j}\right)^{\gamma}\right)\Delta t\rho(0) \\
&= \rho(0)
\end{aligned} \tag{76}
$$

and

$$
\begin{aligned}
g_{i,j}(w) &\leq \left(1 - \Delta t\hat{\theta}_{i,j} - \left(d + k\left(\hat{\theta}_{i,j}\right)^{\gamma}\right)\Delta t\right)\rho\left(w\bar{X}\right) + \Delta t\hat{\theta}_{i,j}\rho\left(w\bar{X}\right) + \left(d + k\left(\hat{\theta}_{i,j}\right)^{\gamma}\right)\Delta t\rho\left(w\bar{X}\right) \\
&= \rho\left(w\bar{X}\right)
\end{aligned} , \tag{77}
$$

proving (47).

Second, from (38) and (42), at each $w > 0$, we obtain

$$
\begin{aligned}
\Phi_{i,j} &= \Phi_{i+1,j} + \Delta t \max_{0 \leq \theta \leq \bar{U}} \left\{ \begin{array}{l} \left(\mathbb{L}_{\theta}\Phi\right)_{i,j} + \theta\bar{h}\bar{W}_{t_i} \\ -\eta\left(\mathbb{L}_{\theta}G\right)_{i,j} + \eta\int_0^{+\infty} p_{W_T}(w)\lambda\left(g_{i+1,j}(w)\right)\left(\mathbb{L}_{\theta}g\right)_{i,j}(w)\mathrm{d}w \end{array} \right\} \\
&\geq \Phi_{i+1,j} + \Delta t\left\{\left(\mathbb{L}_0\Phi\right)_{i,j} - \eta\left(\mathbb{L}_0 G\right)_{i,j} + \eta\int_0^{+\infty} p_{W_T}(w)\lambda\left(g_{i+1,j}(w)\right)\left(\mathbb{L}_0 g\right)_{i,j}(w)\mathrm{d}w\right\} \\
&= \Phi_{i+1,j} + \Delta t\left\{ \begin{array}{l} d\left(\Phi_{i+1,0} - \Phi_{i+1,j}\right) - \eta d\int_0^{+\infty} p_{W_T}(w)\left\{\rho^{(-1)}\left(g_{i+1,0}(w)\right) - \rho^{(-1)}\left(g_{i+1,j}(w)\right)\right\}\mathrm{d}w \\ +\eta\int_0^{+\infty} p_{W_T}(w)\lambda\left(g_{i+1,j}(w)\right) d\left(g_{i+1,0}(w) - g_{i+1,j}(w)\right)\mathrm{d}w \end{array} \right\} \\
&= (1 - d\Delta t)\Phi_{i+1,j} + d\Delta t\Phi_{i+1,0} \\
&\quad +\eta d\Delta t\int_0^{+\infty} p_{W_T}(w)\left\{ \begin{array}{l} \lambda\left(g_{i+1,j}(w)\right)\left(g_{i+1,0}(w) - g_{i+1,j}(w)\right) \\ -\left(\rho^{(-1)}\left(g_{i+1,0}(w)\right) - \rho^{(-1)}\left(g_{i+1,j}(w)\right)\right) \end{array} \right\}\mathrm{d}w \\
&\geq \eta d\Delta t\int_0^{+\infty} p_{W_T}(w)\left\{ \begin{array}{l} \lambda\left(g_{i+1,j}(w)\right)\left(g_{i+1,0}(w) - g_{i+1,j}(w)\right) \\ -\left(\rho^{(-1)}\left(g_{i+1,0}(w)\right) - \rho^{(-1)}\left(g_{i+1,j}(w)\right)\right) \end{array} \right\}\mathrm{d}w \\
&\geq 0
\end{aligned} , \tag{78}
$$

where the second line is due to the that $\theta = 0$ is suboptimal for "max" and we used (44). The last inequality is due to **Lemma C1** which is proven later, and the induction argument.

We now move to the upper bound of the $\Phi_{i,j}$. We set $Z_i = \max_{j=0,1,2,\ldots,N_x} \Phi_{i,j}$. Then, we have

$$
\begin{aligned}
\Phi_{i,j} &= \Phi_{i+1,j} + \Delta t \max_{0\le\theta\le\bar{U}} \left\{ \begin{array}{l} \left(\mathbb{L}_\theta \Phi\right)_{i,j} + \theta \bar{h} \bar{W}_{t_i} \\ -\eta\left(\mathbb{L}_\theta G\right)_{i,j} + \eta \int_0^{+\infty} p_{W_T}(w) \lambda\left(g_{i+1,j}(w)\right)\left(\mathbb{L}_\theta g\right)_{i,j}(w)\mathrm{d}w \end{array} \right\} \\
&= \Phi_{i+1,j} + \Delta t \left\{ \begin{array}{l} \left(\mathbb{L}_\theta \Phi\right)_{i,j} + \theta \bar{h} \bar{W}_{t_i} \\ -\eta\left(\mathbb{L}_\theta G\right)_{i,j} + \eta \int_0^{+\infty} p_{W_T}(w) \lambda\left(g_{i+1,j}(w)\right)\left(\mathbb{L}_\theta g\right)_{i,j}(w)\mathrm{d}w \end{array} \right\}_{\theta=\hat{\theta}_{i,j}} \\
&\le \Phi_{i+1,j} + \Delta t \left\{ \begin{array}{l} \left(\mathbb{L}_\theta \Phi\right)_{i,j} \\ -\eta\left(\mathbb{L}_\theta G\right)_{i,j} + \eta \int_0^{+\infty} p_{W_T}(w) \lambda\left(g_{i+1,j}(w)\right)\left(\mathbb{L}_\theta g\right)_{i,j}(w)\mathrm{d}w \end{array} \right\}_{\theta=\hat{\theta}_{i,j}} \\
&\quad + \bar{U}\bar{h}\Delta t \bar{W}_T
\end{aligned} \quad . \tag{79}
$$

We evaluate each term as follows:

$$
\begin{aligned}
\Phi_{i+1,j} + \Delta t\left(\mathbb{L}_{\hat{\theta}_{i,j}} \Phi\right)_{i,j} &= \Phi_{i+1,j} + \hat{\theta}_{i,j}\Delta t\left(\Phi_{i+1,j-1} - \Phi_{i+1,j}\right) + \left(d + k\left(\hat{\theta}_{i,j}\right)^\gamma\right)\Delta t\left(\Phi_{i+1,0} - \Phi_{i+1,j}\right) \\
&= \left(1 - \left\{\hat{\theta}_{i,j} + d + k\left(\hat{\theta}_{i,j}\right)^\gamma\right\}\Delta t\right)\Phi_{i+1,j} + \hat{\theta}_{i,j}\Delta t \Phi_{i+1,j-1} + \left(d + k\left(\hat{\theta}_{i,j}\right)^\gamma\right)\Delta t \Phi_{i+1,0} \\
&\le \left(1 - \left\{\hat{\theta}_{i,j} + d + k\left(\hat{\theta}_{i,j}\right)^\gamma\right\}\Delta t\right) Z_{i+1} + \hat{\theta}_{i,j}\Delta t Z_j + \left(d + k\left(\hat{\theta}_{i,j}\right)^\gamma\right)\Delta t Z_{i+1} \\
&= Z_{i+1}
\end{aligned} \quad . \tag{80}
$$

We also have

$$
\begin{aligned}
&\Delta t\left\{-\eta\left(\mathbb{L}_\theta G\right)_{i,j} + \eta \int_0^{+\infty} p_{W_T}(w) \lambda\left(g_{i+1,j}(w)\right)\left(\mathbb{L}_\theta g\right)_{i,j}(w)\mathrm{d}w\right\}_{\theta=\hat{\theta}_{i,j}} \\
&= \eta\Delta t\left\{-\left(\mathbb{L}_{\hat{\theta}_{i,j}} G\right)_{i,j} + \int_0^{+\infty} p_{W_T}(w) \lambda\left(g_{i+1,j}(w)\right)\left(\mathbb{L}_{\hat{\theta}_{i,j}} g\right)_{i,j}(w)\mathrm{d}w\right\}
\end{aligned} \quad . \tag{81}
$$

By (47), each term in (81) is evaluated as follows:

$$
\begin{aligned}
-\left(\mathbb{L}_{\hat{\theta}_{i,j}} G\right)_{i,j} &= -\hat{\theta}_{i,j}\int_0^{+\infty} p_{W_T}(w)\left\{\rho^{(-1)}\left(g_{i+1,j-1}(w)\right) - \rho^{(-1)}\left(g_{i+1,j}(w)\right)\right\}\mathrm{d}w \\
&\quad -\left(d + k\left(\hat{\theta}_{i,j}\right)^\gamma\right)\int_0^{+\infty} p_{W_T}(w)\left\{\rho^{(-1)}\left(g_{i+1,0}(w)\right) - \rho^{(-1)}\left(g_{i+1,j}(w)\right)\right\}\mathrm{d}w \\
&\le \hat{\theta}_{i,j}\int_0^{+\infty} p_{W_T}(w)\left\{\rho^{(-1)}\left(\rho\left(w\bar{X}\right)\right) - \rho^{(-1)}\left(\rho(0)\right)\right\}\mathrm{d}w \\
&\quad +\left(d + k\left(\hat{\theta}_{i,j}\right)^\gamma\right)\int_0^{+\infty} p_{W_T}(w)\left\{\rho^{(-1)}\left(\rho\left(w\bar{X}\right)\right) - \rho^{(-1)}\left(\rho(0)\right)\right\}\mathrm{d}w \\
&= \left(\hat{\theta}_{i,j} + d + k\left(\hat{\theta}_{i,j}\right)^\gamma\right)\bar{W}_T\bar{X} \\
&\le \left(\bar{U} + d + k\bar{U}^\gamma\right)\bar{W}_T\bar{X}
\end{aligned} \tag{82}
$$

and similarly,

$$
\begin{aligned}
&\int_0^{+\infty} p_{W_T}(w)\lambda\left(g_{i+1,j}(w)\right)\left(\mathbb{L}_{\hat{\theta}_{i,j}} g\right)_{i,j}(w)\mathrm{d}w \\
&\le \int_0^{+\infty} p_{W_T}(w)\lambda\left(\rho\left(w\bar{X}\right)\right)\left|\left(\mathbb{L}_{\hat{\theta}_{i,j}} g\right)_{i,j}(w)\right|\mathrm{d}w \\
&\le \left(\bar{U} + d + k\bar{U}^\gamma\right)\int_0^{+\infty} p_{W_T}(w)\lambda\left(\rho\left(w\bar{X}\right)\right)\left(\rho\left(w\bar{X}\right) - \rho(0)\right)\mathrm{d}w
\end{aligned} \quad . \tag{83}
$$

Combining (79)-(83) yields

$$\begin{aligned}\Phi_{i,j} \le Z_{i+1} &+ \bar{U}\bar{h}\bar{W}_T \Delta t + \left(\bar{U} + d + k\bar{U}^{\gamma}\right)\eta \bar{W}_T \bar{X} \Delta t \\ &+ \left(\bar{U} + d + k\bar{U}^{\gamma}\right)\eta \Delta t \int_0^{+\infty} p_{W_T}(w)\lambda\left(\rho\left(w\bar{X}\right)\right)\left(\rho\left(w\bar{X}\right) - \rho(0)\right)\mathrm{d}w\end{aligned}. \quad (84)$$

As the right-hand side of (84) is independent of $j$, by a sufficiently large $\bar{\Phi}$ independent of $\Delta t, \Delta x$, we obtain

$$Z_i \le Z_{i+1} + \bar{\Phi}\Delta t \le Z_{N_t} + \bar{\Phi}\left(N_t - i\right)\Delta t = \eta \bar{W}_T \bar{X} + \bar{\Phi}\left(N_t - i\right)\Delta t . \quad (85)$$

This implies that the inequality (46) is satisfied because of the assumption of the induction argument.

□

We finally prove **Lemma C1**.

***Lemma C1***

*For any* $x, y \ge 0$*, strictly increasing and convex (resp., concave)* $\rho : [0, +\infty) \to [0, +\infty)$*, it follows that*

$$\lambda(x)(y - x) - \left(\rho^{(-1)}(y) - \rho^{(-1)}(x)\right) \ge 0 \ \textit{(resp., } \le 0\textit{)}, \quad (86)$$

*where* $\lambda$ *is the derivative of* $\rho^{(-1)}$.

***Proof of Lemma C1***

We only prove the convex case because the concave case is due to a symmetric argument. First, assume that $x, y > 0$. Then, there are unique $w, z > 0$ such that $x = \rho(w)$ and $y = \rho(z)$. Substituting them into the right-hand side of (86) yields

$$\begin{aligned}\lambda(x)(y - x) - \left(\rho^{(-1)}(y) - \rho^{(-1)}(x)\right) &= \frac{\mathrm{d}\rho^{(-1)}(x)}{\mathrm{d}x}\left(\rho(z) - \rho(w)\right) - \left(\rho^{(-1)}(\rho(z)) - \rho^{(-1)}(\rho(w))\right) \\ &= \left.\frac{\mathrm{d}\rho^{(-1)}(x)}{\mathrm{d}x}\right|_{x=\rho(w)}\left(\rho(z) - \rho(w)\right) - (z - w) \\ &= \left(\frac{\mathrm{d}\rho(w)}{\mathrm{d}w}\right)^{-1}\left(\rho(z) - \rho(w)\right) - (z - w) \\ &= \left(\frac{\mathrm{d}\rho(w)}{\mathrm{d}w}\right)^{-1}\left\{\rho(z) - \rho(w) - \frac{\mathrm{d}\rho(w)}{\mathrm{d}w}(z - w)\right\} \\ &\ge 0\end{aligned}, \quad (87)$$

where the last inequality comes from $\frac{\mathrm{d}\rho(w)}{\mathrm{d}w} > 0$ and the convexity of $\rho$. Indeed, if $z > w$, then

$$\rho(z) - \rho(w) - \frac{\mathrm{d}\rho(w)}{\mathrm{d}w}(z - w) = (z - w)\left(\frac{\rho(z) - \rho(w)}{z - w} - \frac{\mathrm{d}\rho(w)}{\mathrm{d}w}\right) > 0 \quad (88)$$

because the tangent line of $\rho$ at $w$ has a slope not larger than $\frac{\rho(z) - \rho(w)}{z - w}$. Similarly, if $z < w$, then

$$\rho(z)-\rho(w)-\frac{\mathrm{d}\rho(w)}{\mathrm{d}w}(z-w)=(w-z)\left(\frac{\mathrm{d}\rho(w)}{\mathrm{d}w}-\frac{\rho(z)-\rho(w)}{z-w}\right)>0\,. \tag{89}$$

Second, assume that $x=0$ and $y>0$. If $\lambda(0)=+\infty$, then we formally obtain (86) due to $y-x>0$ and $\rho^{(-1)}(y)-\rho^{(-1)}(x)$ is bounded; if $\lambda(0)<+\infty$, then the proof of the first case applies. The case $y=0$ and $x>0$ is proven again following the argument in the first case. This completes the proof because $x,y\geq 0$ are arbitrary.

□